\documentclass[a4paper]{article}

\usepackage{setspace} 
\usepackage{subfigure}
\usepackage[margin=0.5in]{geometry}

\usepackage{booktabs} 
\usepackage{array}  
\usepackage{paralist} 
\usepackage{verbatim} 
\usepackage{multicol}[1999/05/25]
\usepackage{graphicx}
\usepackage{epstopdf}
\usepackage{eso-pic} 
\usepackage[toc,page]{appendix}
\usepackage{textpos}
\usepackage{varioref} 
\usepackage{url}
\usepackage{amsmath}
\usepackage{amssymb}
\usepackage{amsfonts}
\usepackage{bm}
\usepackage{float}
\usepackage{doi}
\newcommand{\e}[1]{{\mathbb E}\left[ #1 \right]}

\def\XXint#1#2#3{{\setbox0=\hbox{$#1{#2#3}{\int}$}
		\vcenter{\hbox{$#2#3$}}\kern-.5\wd0}}

\begin{document}

\title{Stochastic analysis of a full system of two competing populations in a chemostat}

\author{
	Dimitrios Voulgarelis$^{\ast}$\thanks{$^\ast$Corresponding author. Email: dimitrios.voulgarelis.14@ucl.ac.uk}, Ajoy Velayudhan and Frank Smith
}
\maketitle

	\begin{abstract}
	This paper formulates two 3D stochastic differential equations (SDEs) of two microbial populations in a chemostat competing over a single substrate. The two models have two distinct noise sources. One is general noise whereas the other is dilution rate induced noise. Nonlinear Monod growth rates are assumed and the paper is mainly focused on the parameter values where coexistence is present deterministically. Nondimensionalising the equations around the point of intersection of the two growth rates leads to a large parameter which is the nondimensional substrate feed. This in turn is used to perform an asymptotic analysis leading to a reduced 2D system of equations describing the dynamics of the populations on and close to a line of steady states retrieved from the deterministic stability analysis. That reduced system allows the formulation of a spatially 2D Fokker-Planck equation which when solved numerically admits results similar to those from simulation of the SDEs. Contrary to previous suggestions, one particular population becomes dominant at large times. Finally, we briefly explore the case where death rates are added.
	\end{abstract}

	\section{Introduction}
	
	\subsection{Motivation}
	
	For decades the growth of bacterial/cell populations has been a subject of great interest to modellers. The reason behind the popularity of these systems is of course the industrial and ecological importance of competing population growth processes as well as the richness and complexity of the dynamics arising from even simple systems of a few competing organism. When exploring these systems, coexistence of the different populations is of great significance. One of the best known papers on analysis of coexistence of competing populations is the paper by Stephanopoulos et al. \cite{stephanopoulos_stochastic_1979} where they explored the dynamics of two microbial populations competing for a single substrate. Before that paper it was proven in \cite{aris_dynamics_1977} that if the substrate concentration is kept at the break-even point of the two populations and the dilution rate constant at that value then both populations can coexist. The result was also generalized and proven for multiple competing populations in \cite{hsu1978limiting} with the use of Lyapunov functions. However that situation presents a knife-edge event where the extinction of one or of the other population occurs if the dilution rate diverges from this exact value \cite{el_hajji_practical_2009}. The main result by Stephanopoulos et al. was that this extinction occurs due to the noise present in the control of the dilution rate in every chemostat. The interesting conclusion was that either population can become extinct depending on the value of the dilution rate and the initial conditions. \\
	
	A chemostat is an automated bioreactor in which spent medium which contains metabolic products, microorganisms and left over nutrients is continuously removed while fresh medium is added at the same rate to keep the volume constant \cite{novick1950description}. That rate is called the dilution rate and in the case where it is smaller than the growth rate of the micro-organism that micro-organism will grows. The chemostat provides a powerful means of systematically investigating how growth rate impacts processes of the cells such as gene expression and metabolism and the regulatory networks that control the rate of cell growth. Moreover, cells grown in chemostat for generations can be used to study their adaptive evolution in environmental conditions that limit cell growth \cite{ziv_use_2013}. One of the most important characteristics of the chemostat for multiple microorganism populations competing over a single substrate is the Competitive Exclusion Principle (CEP). Per the CEP in the above scenario only one population will survive. More specifically the one that has the lowest break-even concentration will survive while the other will be led to extinction. The break-even concentration is the concentration of the nutrient such that the specific growth of a microorganism is equal to the dilution rate. A great number of papers have focused on proving the CEP for different growth function assumptions and removal rates. Most of the papers use deterministic equations to describe the evolution of populations in the chemostat while a few have recently addressed what happens when stochasticity is taken into account with either linear growth rates \cite{xu_competition_2016} or with only a single population \cite{ji_dynamics_2014,xu2015analogue}.
	
	\subsection{Aim}
	
	The aim of this paper is to explore the idea of coexistence by simulating the dynamics of the full equations for two microbial populations and one substrate for non-linear Monod growth functions with general noise as well as dilution rate induced noise as explored in \cite{stephanopoulos_stochastic_1979}. The rest of the paper is divided into \textbf{materials and methods}, where we present the stochastic version of the full model, for both cases, in the form of a set of three stochastic differential equations (SDEs) of the Langevin type. Beforehand an asymptotic analysis, which is performed for the case where the substrate feed is large to aid our understanding of the system, shows an intricate structure within the dynamics and provides a simplified two-dimensional version from which we can derive and numerically solve the Fokker-Planck equation readily. Finally we examine the case were death rates are added to the model solely for the dilution rate noise case. The next section after that is \textbf{Results and discussion}. Here, the equations are numerically solved and simulated for the parameter values of the same two microbial populations used in \cite{stephanopoulos_stochastic_1979}. Following the results section, our work and findings are summarized in \textbf{Further discussion} and finally in the last section named \textbf{Conclusion} we raise possible issues as well as possible extensions to our work.
	
	\section{Materials and Methods}
	
	For two populations in a chemostat competing over a single substrate the dimensionless equations are given by:
	
	\begin{align}
	\frac{dx}{dt}&=x(f(z)-\theta),\\
	\frac{dy}{dt}&=y(g(z)-\theta),\\
	\frac{dz}{dt}&=\theta(z_f-z)-xf(z)-yg(z).
	\end{align}
	
	\noindent Here $f(z), g(z)$ represent the dimensionless growth rates given by the following equations:
	
	\begin{align}
	f(z)=\frac{a_1z}{b_1+z},\\
	g(z)=\frac{a_2z}{b_2+z}.
	\end{align}
	
	\noindent A list of the parameters and their definitions is given in Table 1.\\
	
	\begin{table}[ht]
		\centering
		\caption{Parameters and their values for E. coli and Spirillum sp. respectively.}
		{\begin{tabular}[l]{@{}lcccccc}\toprule
				Parameter (dimensionless) & Definition & Value & Reference \\ \hline
				$\theta$ & dilution rate & varying & -\\ 
				$z_f$ & substrate feed & 15000 & -\\ 
				$a_i$ & maximum growth rate & 2.911, 1.636 & \cite{stephanopoulos_stochastic_1979}\\ 
				$b_i$ & Michaelis constant & 1.911, 0.636 & \cite{stephanopoulos_stochastic_1979}\\ 
				$\sigma$ & noise intensity & varying & our model\\ 
				$\theta_0$ & dilution rate mean & 1 & -\\ \hline
		\end{tabular}}
		\label{Table2}
	\end{table}
	
	The non-dimensionalisation was performed around the break-even concentration of the substrate assuming that there is such. In the case of the parameter values used in \cite{stephanopoulos_stochastic_1979} which will also be used here, there is such a point. In the dimensionless system the two growth rates break even when $z=1$ in which case $f(z),g(z)$ are also equal to one. In order to have coexistence of the two populations the value of the dimensionless dilution rate, $\theta$, must be one. Then it can be shown using linear stability that there is a line of steady states given by $y=z_f-x-1$ \cite{stephanopoulos_stochastic_1979}. \\
	
	Equations (1-3) were simplified in \cite{stephanopoulos_stochastic_1979} and the system reduced to one dimension before introducing the noise term in the dilution rate. In our analysis, by contrast, the noise term is introduced in the full equations without further simplifications, first, and computational studies are made; afterwards a self-consistent asymptotic treatment is also applied to complement the numerical approach and provide further comparisons. \\
	
	\subsection{Stochastic Langevin equations}

	It was been show in \cite{imhof2005exclusion} that in a chemostat system of the form (1-3) stochastic effects can be added as follows:
	
	\begin{align}
	dx&=x(f(z)-\theta_0)dt + \sigma_1 x dW_1(t) ,\\
	dy&=y(g(z)-\theta_0)dt + \sigma_2 y dW_2(t),\\
	dz&=\big[\theta_0(z_f-z)-xf(z)-yg(z)\big]dt + \sigma_3 z dW_3(t).
	\end{align}
	\\
	
	\noindent Here $W_i$ are independent Wiener processes (Brownian motions).\\
	
	In the case of stochasticity being solely due to random fluctuation in the dilution rate the equations are different. Here, the dilution rate $\theta$ fluctuates around a mean value and so:
	
	\begin{equation}
	\theta=\theta_0+\zeta(\tau),
	\end{equation}
	\\
	
	\noindent Here, $\zeta(\tau)$ is a Gaussian random noise. Substituting that back into (1-3), a system of stochastic differential equations is found:
	
	\begin{align}
	dx&=x(f(z)-\theta_0)dt - \sigma x dW(t) ,\\
	dy&=y(g(z)-\theta_0)dt - \sigma y dW(t),\\
	dz&=\big[\theta_0(z_f-z)-xf(z)-yg(z)\big]dt + \sigma(z_f-z) dW(t).
	\end{align}
	\\
	
	\noindent Here $W$ is a Wiener process and the formal derivative of the Gaussian noise. The Wiener process increments are the same in each of the equations because there is a unique source of noise, i.e. the variation in the dilution rate.

	\subsection{Asymptotic analysis for large $z_f$}
	
	According to \cite{stephanopoulos_stochastic_1979} the Michaelis constants for most populations are of the order of a few milligrams per litre, and as a result the break-even value of the substrate concentration at which the two specific growth rate curves intersect is of
	similar order. The concentration of the substrate in the feed is of the order of several grams per litre, so that $z_f$ is of the order of several thousands. That gives a natural large parameter in the system which we can use to perform an asymptotic analysis. The large value of $z_f$ is of great importance to the analysis of \cite{stephanopoulos_stochastic_1979} since it is used to suggest that the movement along the line of steady states is slow whereas the movement of the system from any point close to the line towards the line itself is very fast. Hence, according to this suggestion the behaviour of the system around the line can be ignored and the system modelled on the line. The reason behind the present analysis is to understand the system better close to the line of steady states and determine whether we can justify the reduction of the system to one dimension or not. In fact using the results from the asymptotic analysis indicates we can only reduce the system to two dimensions near the line of steady states and we then numerically solve the corresponding Fokker-Planck PDE readily.\\
	
	We can investigate what happens for $x(0), y(0)$ of O(1) and larger: all such cases are of some concern and validity, and they pass through a number of successive temporal stages. A central case however is found to occur for specific combinations of the initial conditions on x, y, namely
	
	\begin{equation}
	x(0)=M_1, y(0)=M_2.
	\end{equation}
	
	\noindent The initial values $M_1,M_2$ are considered later.
	
	\subsubsection{Stage 1}
	Since all variables start at order less than $z_f$ we can see from the system of equations that in order to balance the equation we need to introduce a fast time $T=z_ft$. This will give us $\frac{d}{dt} =z_f*\frac{d}{dT}$. Hence, expanding the variables as:
	
	\begin{align*}
	x &=  x_0 + \epsilon x_1 + \ldots,\\
	y &= y_0 +  \epsilon y_1 + \ldots,\\
	z &= z_0 + \epsilon z_1 + \ldots,
	\end{align*}
	
	\noindent and taking the O($z_f$) terms gives:
	
	\begin{align*}
	\frac{dx_0}{dT} &= 0,\\
	\frac{dY_0}{dT} &= 0,\\
	\frac{dz_0}{dT} & = \theta.
	\end{align*}
	
	The above means that in the first stage x,y will remain constant while z will increase linearly with time, $z \approx z_0 + \theta z_ft$ for a period of O(1/$z_f$).\\
	
	\subsubsection{Stage 2}
	
	Stage 2 has t of O(1). After the end of stage 1, $z$ has been increased to O($z_f$) whereas x and y remained constant. In stage 2, x, y grow exponentially at order-unity rates while z saturates near the value $z_f$ at large t. Now we can introduce a new expansion for $z$:
	
	\begin{equation*}
	z= z_f\bar{z} + \ldots,
	\end{equation*}
	
	Furthermore we have constants $b_i \sim O(1)$ and hence $b_i << z$ which means that $f(z) \sim a_1$ and $g(z) \sim a_2$. Using these we have to leading order:
	
	\begin{align*}
	\dot{x} &= x(a_1-\theta),\\
	\dot{y} &= y(a_2-\theta),\\
	\dot{\bar{z}} &= \theta(1-\bar{z}),
	\end{align*}
	
	\noindent Thus the solutions are
	
	\begin{align}
	x &= x(0)e^{(a_1-\theta)t},\\
	y &= y(0)e^{(a_2-\theta)t},\\
	z &= z_f(1-C_0e^{-\theta t}),
	\end{align}
	
	\begin{figure}[H]
		\begin{center}
			\resizebox*{14cm}{!}{\includegraphics{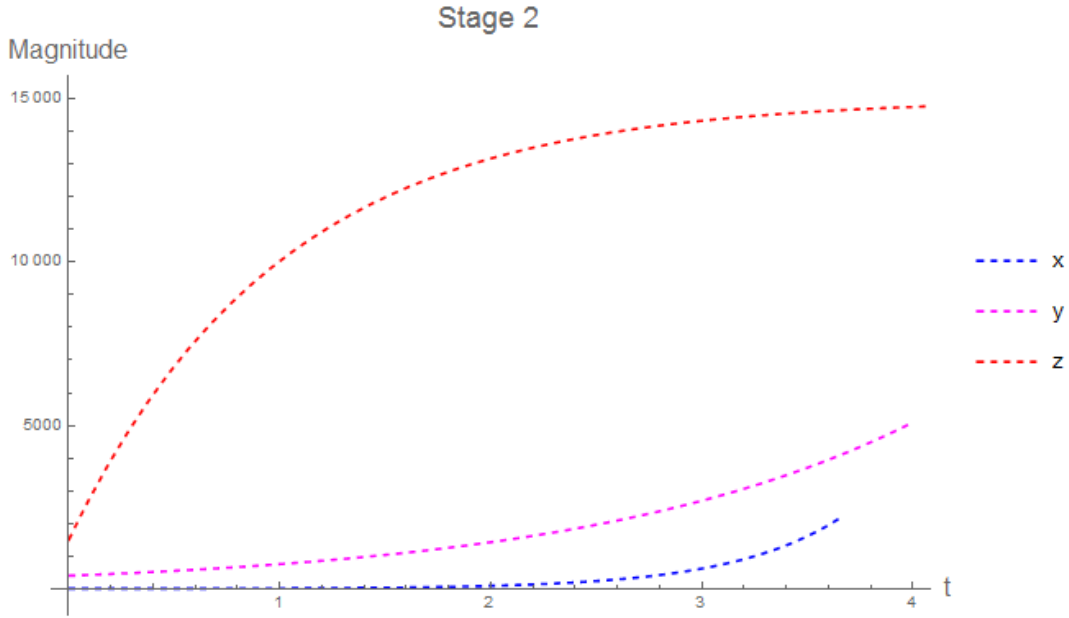}}
			\caption{Plot of $x,y,z$ for the second stage.} \label{fig3.6}
		\end{center}
	\end{figure}
	
	Since both x and y increase exponentially they will become large quickly and start affecting $z$. We want as a central case both x and y to have an almost equal contribution to the dynamics so we need to pick the appropriate initial conditions for y, i.e. y(0). Namely, what we examine is $x,y \sim z_f$ at the same time. 
	
	\begin{equation}
	M_1e^{(a_1-\theta)t}\sim z_f, M_2e^{(a_2-\theta)t}\sim z_f.
	\end{equation}
	
	\subsubsection{Stage 3}
	
	In this stage x, y interact directly with z. This interactive stage has the form:
	
	\begin{equation}
	t=L+t',(x,y,z)=z_f(x',y',z')+\ldots,
	\end{equation}
	
	\noindent where the primed quantities are of order unity and the large parameter $L$ satisfies
	
	\begin{equation}
	L=(a_1-\theta)^{-1}\log{(z_f/M_1)}=(a_2-\theta)^{-1}\log{(z_f/M_2)},
	\end{equation}
	
	\noindent owing to (17).\\
	
	Given the expansion (18), the system (1)-(3) reduces to the linear equations
	
	\begin{align}
	\frac{dx'}{dt}&=x'(a_1-\theta),\\
	\frac{dy'}{dt}&=y'(a_2-\theta),\\
	\frac{dz'}{dt}&=\theta(1-z')-x'a_1-y'a_2,
	\end{align}
	
	\noindent in stage 3. The solution for $x',y',z'$ here yields, after matching with the previous stages,
	
	\begin{align}
	x'&= \mu_1e^{(a_1-\theta)t'},\\
	y'&= \mu_2e^{(a_2-\theta)t'},\\
	z'&=1-C_3e^{-\theta t'}-\mu_1e^{(a_1-\theta)t'}-\mu_2e^{(a_2-\theta)t'},
	\end{align} 
	
	\begin{figure}[H]
		\begin{center}
			\resizebox*{14cm}{!}{\includegraphics{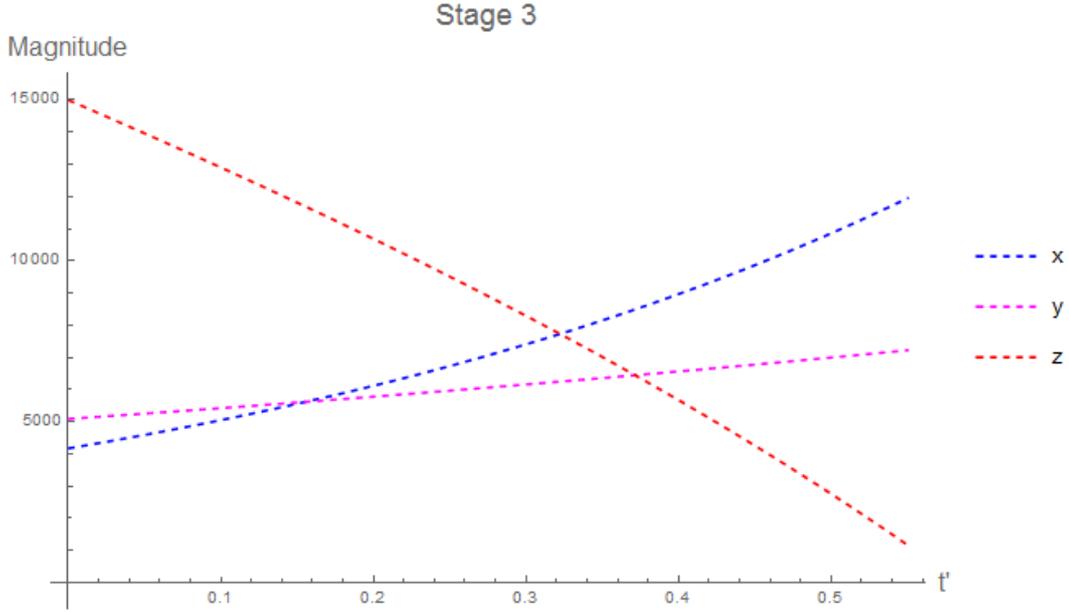}}
			\caption{Plot of $x,y,z$ for the third stage.} \label{fig3.6}
		\end{center}
	\end{figure}
	
	\noindent with the coefficients $\mu_1, \mu_2$ of the two growing exponentials being notable. From matching with the previous stage we deduce that $C_3=0$. As $t′$ increases from $-\infty$ the solution $z'$ increases monotonically at first but then at a finite time $z'$ achieves a maximum value, after which $z'$ decreases monotonically and reaches zero within a further finite time, say at $t' = t_0'$.\\
	
	To find the coefficients $\mu_1, \mu_2$ as a function of $M_1,M_2$ we require the following:
	
	\begin{align}
	M_1e^{(a1-\theta)L}&=\mu_1z_f,\\
	M_2e^{(a2-\theta)L}&=\mu_2z_f,
	\end{align}
	
	\noindent which by fixing $M_1,M_2$ and finding $L$ can admit the coefficients for stage3.\\
	
	A significant condition for $z'$ to reach that zero comes straight from (22), which implies that
	
	\begin{equation}
	\theta-x_0'a_1-y_0'a_2 < 0.
	\end{equation}
	
	\noindent based on the assumption that $z'$ approaches zero from above at $t' = t_0'-$, i.e. $dz'/dt'$ must then be negative. The quantities $x_0', y_0'$ in (28) are the values of $x', y'$ at $t' = t_0'–$.
	
	\subsubsection{Stage 4}
	
	Stage 4 arises when z reduces to the order of unity. This is a rapid decrease stage in which we have
	
	\begin{equation}
	t = L + t_0' + z_f^{-1}T, (x, y, z) = (z_f x_0', z_f x_0', Z) +\ldots,
	\end{equation}
	
	Thus $x, y$ remain constant to leading order. The governing system (1)-(3) produces evolution equations for the perturbations in $x, y$, while for $Z$ we find that
	
	\begin{equation}
	\frac{dZ}{dt}=\theta-x_0'f(Z)-y_0'g(Z),
	\end{equation}
	
	\noindent Here $f(Z), g(Z)$ are non-trivial, being respectively $a_1 Z / ( b_1 + Z)$, $a_2 Z / ( b_2 + Z)$. Matching with stage 3 at large negative T we again find that condition (28) needs to be satisfied in order for Z to decrease at the start of stage 4. That leads to an evolution of $Z$ towards the constant value $Z_{\infty}$ at large positive T. The constant satisfies:
	
	\begin{equation}
	x_0'f(Z_{\infty})+y_0'g(Z_{\infty})=\theta,
	\end{equation}
	
	\noindent due to (30).\\
	
	\begin{figure}[H]
		\begin{center}
			\resizebox*{14cm}{!}{\includegraphics{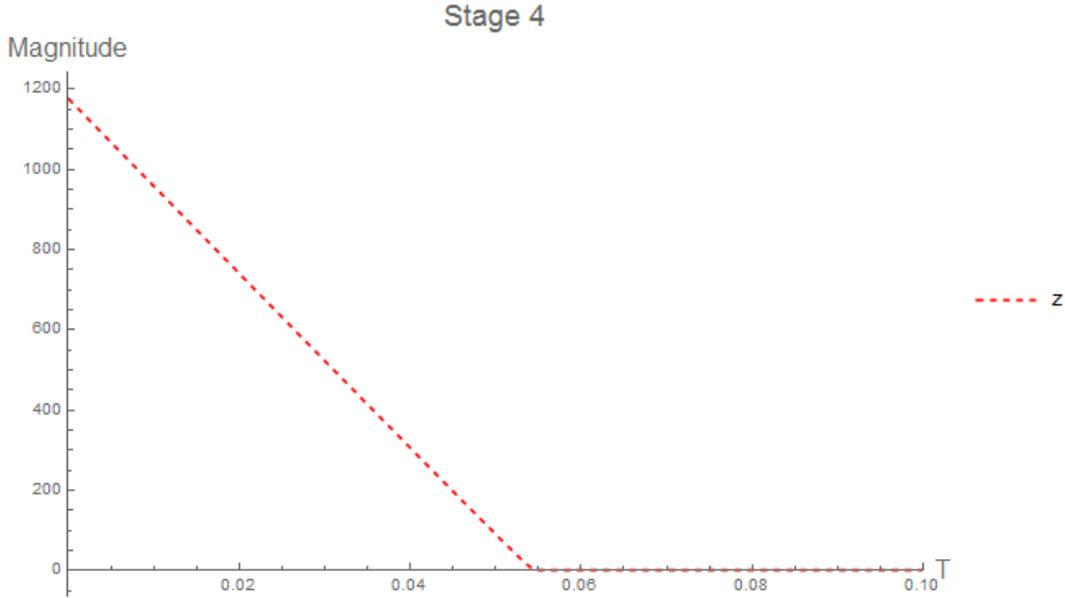}}
			\caption{Plot of $z$ for the fourth stage 4.} \label{fig3.6}
		\end{center}
	\end{figure}

	We omitted the $x,y$ plot as it is mostly trivial especially in the case of $\theta=1$ where they reach the steady state by the end of stage 4.	 
	\subsubsection{Stage 5}
	
	Stage 5 is the final stage. Here x, y are of O($z_f$), whereas z is O(1) and the typical time variation is of order unity. Hence
	
	\begin{equation}
	t = L + t_0' + \bar{t}, (x, y, z) = (z_f \bar{x}, z_f \bar{y}, \bar{z}) + \ldots,
	\end{equation}
	
	The equations in this stage reduce to
	
	\begin{align}
	\dot{\bar{x}} &= \bar{x}(f(\bar{z})-\theta),\\
	\dot{\bar{y}} &= \bar{y}(g(\bar{z})-\theta),\\
	\theta &= \bar{x}f(\bar{z}) +  \bar{y}g(\bar{z}),
	\end{align}
	
	\begin{figure}[H]
		\begin{center}
			\resizebox*{14cm}{!}{\includegraphics{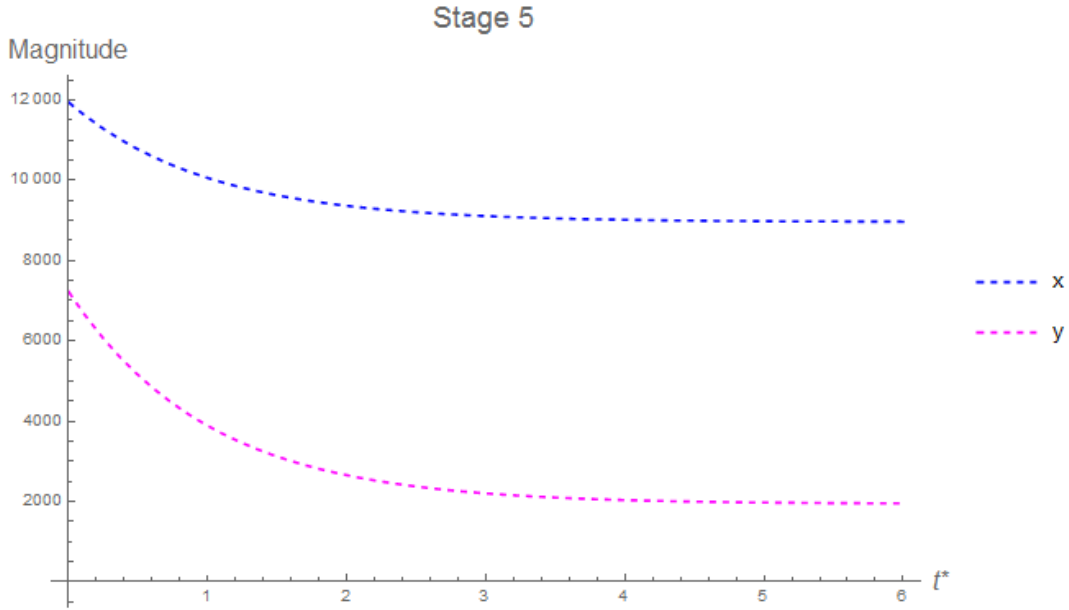}}
			\caption{Plot of $x,y$ for the fifth stage.} \label{fig3.6}
		\end{center}
	\end{figure}
	
	\noindent where the feedback effects $f(\bar{z}) = a_1 \bar{z} / ( b_1 + \bar{z}), g(\bar{z}) = a_2 \bar{z} / ( b_2 + \bar{z})$ are nontrivial again. The initial conditions are:
	
	\begin{equation}
	(\bar{x}(0), \bar{y}(0), \bar{z}(0)) = (x_0' , y_0', Z_{\infty}),
	\end{equation}
	
	\noindent from matching with the previous stage. We observe that the initial conditions satisfy the restriction
	
	\begin{equation}
	\bar{x}(0) a_1 + \bar{y}(0) a_2 > \theta ,
	\end{equation}
	
	\noindent in view of the requirement (28) in an earlier stage.\\ 
	
	In the differential-algebraic equations (DAEs) above $\bar{z}$ is given by solving the algebraic expression whereas x and y are given by the regular differential equations of the full system. For the following analysis as well as the Fokker-Planck equation we will use the variables $\bar{x},\bar{y}$ and $\bar{z}$.The equation for $\bar{z}$ is quadratic so we end up with two possible solutions for $\bar{z}$. If we make a contour plot of the solution for $\bar{x},\bar{y}$ varying between 0 and 1 we can see that one solution yields only negative results so it can be dropped. Keeping the other solution though yields a line of singularities given by:
	
	\begin{equation}
	\bar{y} = \frac{\theta}{a_2} - \frac{a_1}{a_2}\bar{x},
	\end{equation}
	
	\begin{figure}[ht!]
		\begin{center}
			\subfigure[Quadratic solution 1.]{
				\resizebox*{6cm}{!}{\includegraphics{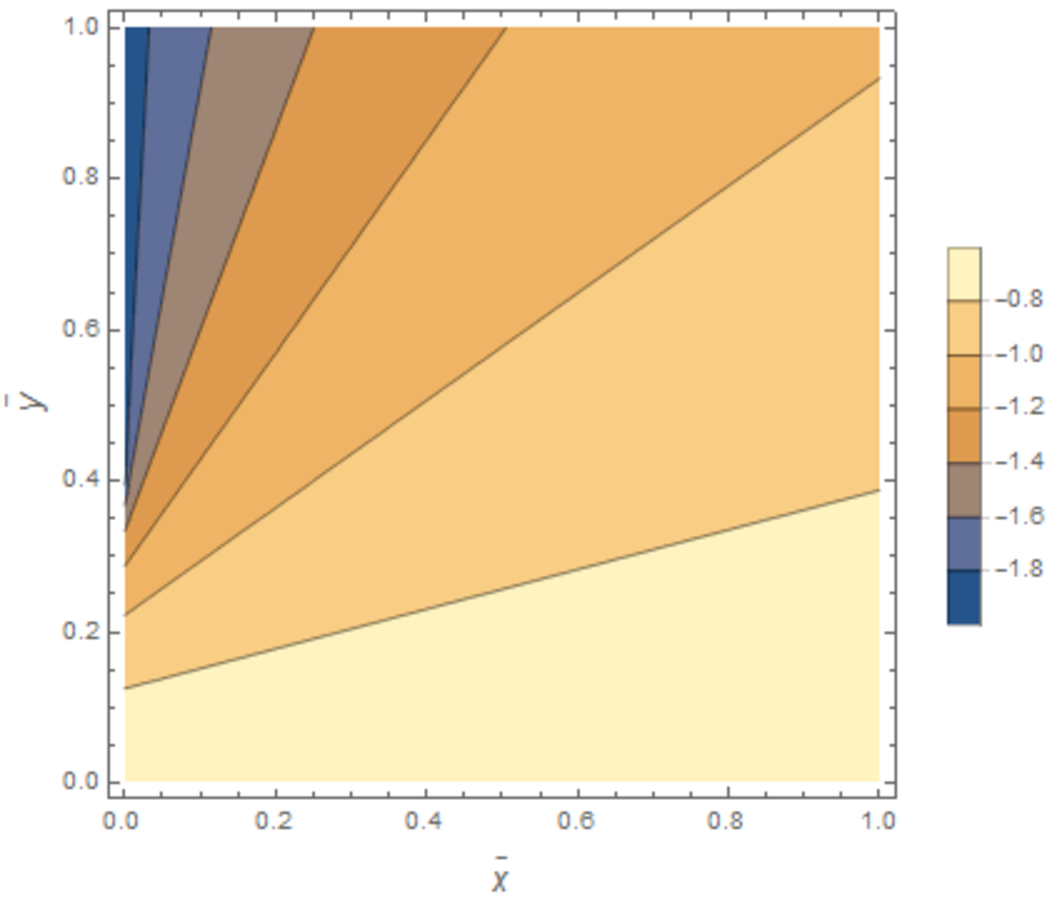}}}\hspace{5pt}
			\subfigure[Quadratic solution 2.]{
				\resizebox*{6cm}{!}{\includegraphics{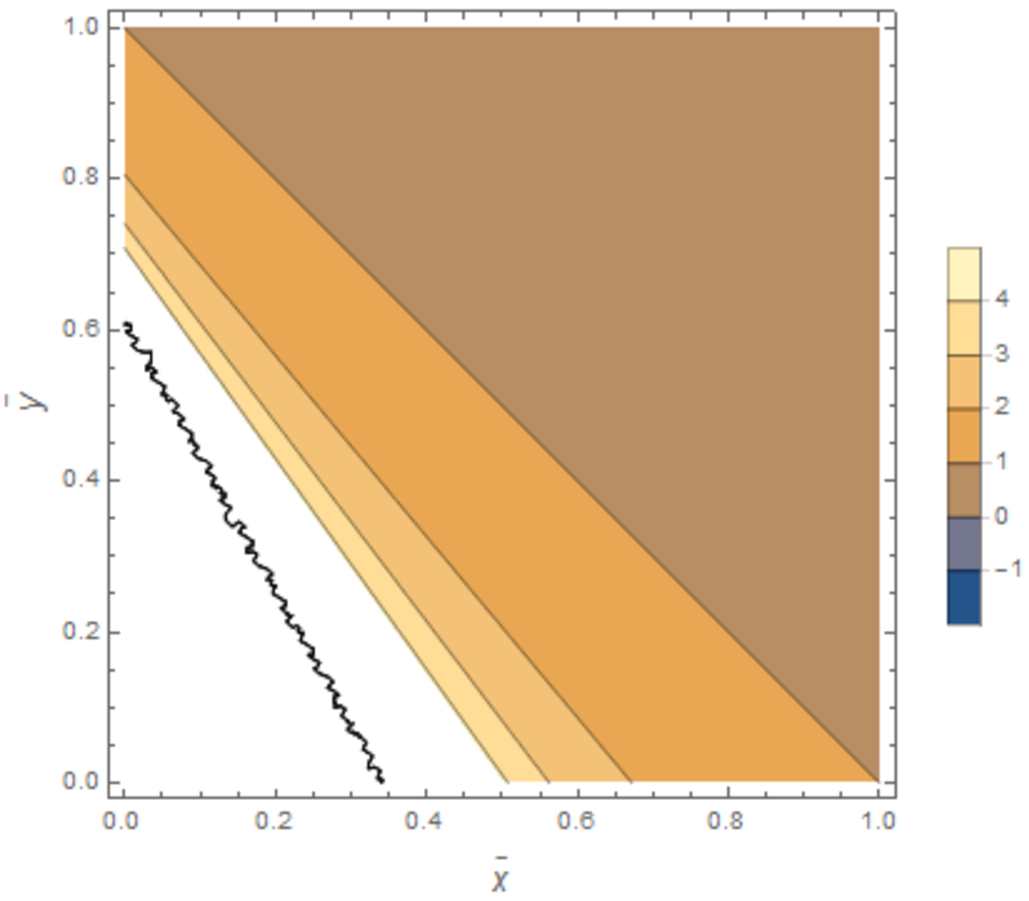}}}
			\caption{Contour plot of the two solutions for z from the quadratic equation derived in stage 4.} \label{figure10}
		\end{center}
	\end{figure}
	
	\noindent Hence, we need to solve the system from this line onwards. To do that we first need to make sure that in stage 5 the system is already far from that line. This can be easily deduced from the fact that in stage 5 the system already starts away from the line of singularities due to (37).\\
	
	The asymptotic approach provides a firm basis for the deterministic case but it also points towards use of the apparently self-consistent approximations in the present work in terms of increased understanding of the stochastic differential and Fokker-Planck equations. As far as we know no such basis exists as yet for the approximations in \cite{stephanopoulos_stochastic_1979}.
	
	\subsection{Fokker-Planck equation of stage 5 Langevin equations with the assumption of constant substrate source}
	
	Deriving the Langevin equations for stage 5 so as to find the respective Fokker-Planck equations is trivial and identical to deriving the stochastic equations for the full system. For the Fokker-Planck equations we use the methodology shown in Appendix A and we obtain two different equations depending on the noise assumptions. More specifically in the case where the Brownian motion increments are the same we gain an extra term. In the process of deriving the Fokker-Planck equations we make a significant assumption/simplification which is that the algebraic equation for z has no noise. The only variation of $\bar{z}$ due to noise comes from the variability of $\bar{x},\bar{y}$. Without that assumption we cannot readily derive a compact Fokker-Planck equation for stage 5. We will show that this assumption makes z vary significantly less than how much it would in the case of noise in the algebraic equation. This doesn't seem to affect the qualitative results and it provides a compromise between no z variation in \cite{stephanopoulos_stochastic_1979} and full $\bar{z}$ stochastic variation since from simulation of the stochastic equations for the full system we will see that the smaller the noise in the z equation the more noisy $x,y$ are close to the line. As a result the system takes more time to reach a steady-state. Hence, the limit is to consider the noise intensity of equation (8 \& 12) as zero and recover equation (3). Which when close to the line of steady states will in turn be replaced by (35) due to the asymptotic analysis \\
	
	The biological premise of not including noise effects on the equation from $\bar{z}$ arises if we consider a chemostat where a constant source of substrate is provided equal to $\theta_0 z_f$. In that case we can suppose that the change in the dilution rate only affects the substrate evolution indirectly through the the stochastic effect on $x,y$ populations.\\
	
	Using the general formula for the derivation of the Fokker-Planck equation from Appendix A, the Fokker-Planck equations for the simplified 2D system of the fifth stage, for the general noise and dilution rate noise respectively, are found to be:
	
	\begin{align}
	\begin{split}
	\frac{\partial}{\partial t}p(\bar{x},\bar{y},t|\bar{x}_0,\bar{y}_0,t_0)=&-\frac{\partial}{\partial \bar{x}}\big[\bar{x}(f(\bar{x},\bar{y})-\theta_0)p(\bar{x},\bar{y},t|\bar{x}_0,\bar{y}_0,t_0)\big]\\
	&-\frac{\partial}{\partial \bar{y}}\big[\bar{y}(g(\bar{x},\bar{y})-\theta_0)p(\bar{x},\bar{y},t|\bar{x}_0,\bar{y}_0,t_0)\big]\\
	&+\frac{\sigma_1^2}{2}\frac{\partial^2}{\partial \bar{x}^2}\big[\bar{x}^2p(\bar{x},\bar{y},t|\bar{x}_0,\bar{y}_0,t_0)\big]+\frac{\sigma_2^2}{2}\frac{\partial^2}{\partial \bar{y}^2}\big[\bar{y}^2p(\bar{x},\bar{y},t|\bar{x}_0,\bar{y}_0,t_0)\big],
	\end{split}
	\end{align}
	
	\begin{align}
	\begin{split}
	\frac{\partial}{\partial t}p(\bar{x},\bar{y},t|\bar{x}_0,\bar{y}_0,t_0)=&-\frac{\partial}{\partial \bar{x}}\big[\bar{x}(f(\bar{x},\bar{y})-\theta_0)p(\bar{x},\bar{y},t|\bar{x}_0,\bar{y}_0,t_0)\big]\\
	&-\frac{\partial}{\partial \bar{y}}\big[\bar{y}(g(\bar{x},\bar{y})-\theta_0)p(\bar{x},\bar{y},t|\bar{x}_0,\bar{y}_0,t_0)\big]\\
	&+\frac{\sigma^2}{2}\frac{\partial^2}{\partial \bar{x}^2}\big[\bar{x}^2p(\bar{x},\bar{y},t|\bar{x}_0,\bar{y}_0,t_0)\big]+\frac{\sigma^2}{2}\frac{\partial^2}{\partial \bar{y}^2}\big[\bar{y}^2p(\bar{x},\bar{y},t|\bar{x}_0,\bar{y}_0,t_0)\big]\\
	&+\sigma^2\frac{\partial^2}{\partial \bar{x}\partial \bar{y}}\big[\bar{x}\bar{y}p(\bar{x},\bar{y},t|\bar{x}_0,\bar{y}_0,t_0)\big].
	\end{split}
	\end{align}
	
	\noindent Here $p(\bar{x},\bar{y},t|\bar{x}_0,\bar{y}_0,t_0)$ is the probability density of the system being in state $\bar{x},\bar{y}$ at time t assuming that it started in state $\bar{x}_0,\bar{y}_0$ at time $t_0$. Moreover, $\sigma$ is the noise intensity.
	
	\subsubsection{Numerical solution formulation, boundary and initial conditions}
	
	To simulate low noise intensity a very fine grid mesh was used in the polygon areas created by cutting the corner below the line of singularities mentioned in stage 5 in \textbf{2.2.5} and extending the area as far $\bar{x}=3$ and $\bar{y}=3$ which would correspond to $x=y=3z_f$. The corner cut is the area below the line defined as:
	
	\begin{equation}
	y = 1/a_2 - a_1 x/a_2 + 10^{-2}
	\end{equation}
	
	\noindent That is to avoid the numerical errors occurring by being too close to the line of singularities.\\
	
	The mesh used is the following:
	
	\begin{figure}[H]
		\begin{center}
			\resizebox*{10cm}{!}{\includegraphics{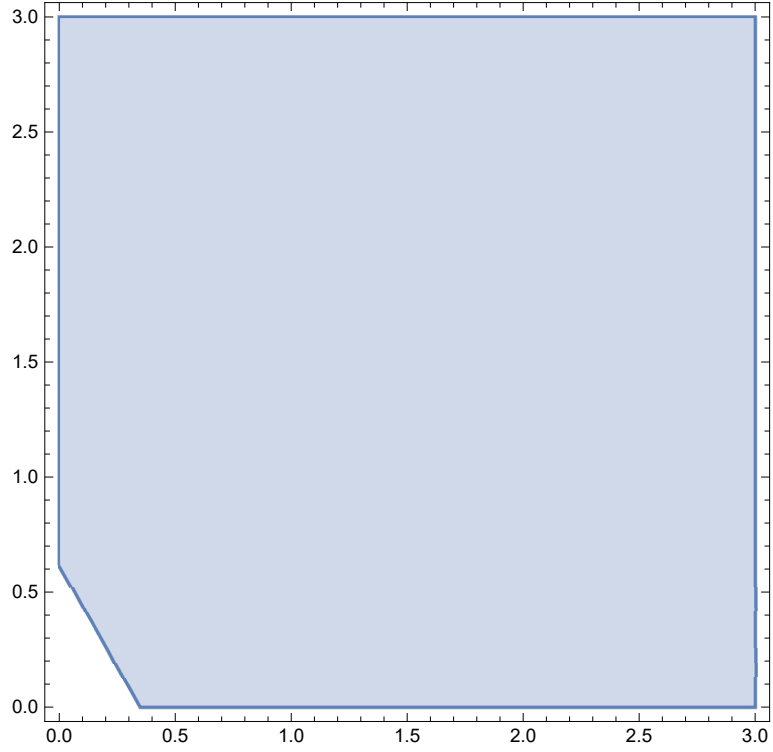}}
			\caption{2D mesh for the solution of the Fokker-Planck equation.} 
		\end{center}
	\end{figure}
	
	The initial condition for the numerical solution of the Fokker-Planck equation is a two dimensional Gaussian with means $\mu_1=\mu_2=0.5$ and standard deviations of $\sigma_1=\sigma_2=0.05$:
	
	\begin{equation*}
	p(x,y,0)=\frac{1}{2\pi\sigma_1\sigma_2}e^{-\frac{1}{2}(\frac{(x-\mu_1)^2}{\sigma_1^2}+\frac{(y-\mu_2)^2}{\sigma_2^2})},
	\end{equation*}
	
	The boundary conditions are such that the probability density is zero at the $\bar{x}=3,\bar{y}=3$ boundary and there is no flux out of the region so that the total probability remains equal to 1. For the boundary at $x=3,y=3$ we impose $p(3,y,0)=p(x,3,0)=0$ as well as $\frac{\partial p}{\partial x}|_{x=3}=0,\frac{\partial p}{\partial y}|_{y=3}=0$. For the $\bar{x}=0$ or $\bar{y}=0$ boundary we do not need to specify $p=0$ boundary conditions as the form of the Fokker-Planck equation at the boundary implies that if $p$ is initially zero there it will remain so. To see that we need only write down (see below) the Fokker-Planck at the $y=0$ boundary for one of the noise cases as the other will admit something similar.

	\begin{align*}
	\frac{\partial}{\partial t}p(\bar{x},0,t|x_0,\bar{y}_0,t_0)=&-\frac{\partial}{\partial \bar{x}}\big[\bar{x}(f(\bar{x},0)-\theta_0)p(\bar{x},0,t|\bar{x}_0,\bar{y}_0,t_0)\big]\\
	&-(g(\bar{x},0)-\theta_0)p(\bar{x},0,t|\bar{x}_0,\bar{y}_0,t_0)\\
	&+\frac{\sigma^2}{2}\frac{\partial^2}{\partial \bar{x}^2}\big[\bar{x}^2p(\bar{x},0,t|\bar{x}_0,\bar{y}_0,t_0)\big]+0\sigma^2p(\bar{x},0,t|\bar{x}_0,\bar{y}_0,t_0)\\
	&+\sigma^2\frac{\partial}{\partial \bar{x}}\big[\bar{x}p(\bar{x},0,t|\bar{x}_0,\bar{y}_0,t_0)\big],
	\end{align*}
	
	\noindent We can see that there is no contribution to the evolution of $p$ on that boundary from values of $y \neq 0$ and that if $p=0$ initially it remains so as $\frac{\partial p}{\partial t}=0$.\\
	
	Finally we require that $p$ is zero on the diagonal line (38) and that the flux at this line is also zero. To enforce that boundary condition we first need to calculate the probability flux. This is given by:
	
	\begin{equation}
	\partial_t p + \nabla\cdot\bf{J}=0,
	\end{equation}
	
	\noindent which from (39,40) admits:
	
	\begin{align}
	\begin{split}
	\bf{J}&=\bigg(\bar{x}(f(\bar{x},\bar{y})-\theta_0)p-\frac{\sigma^2}{2}\frac{\partial}{\partial \bar{x}}\big[\bar{x}^2p\big],\\
	&\bar{y}(g(\bar{x},\bar{y})-\theta_0)p-\frac{\sigma^2}{2}\frac{\partial}{\partial \bar{y}}\big[\bar{y}^2p\big]\bigg),
	\end{split}
	\end{align}
	
	\begin{align}
	\begin{split}
	\bf{J}&=\bigg(\bar{x}(f(\bar{x},\bar{y})-\theta_0)p-\frac{\sigma^2}{2}\frac{\partial}{\partial \bar{x}}\big[\bar{x}^2p\big]-\sigma^2\frac{\partial}{\partial \bar{y}}\big[\bar{x}\bar{y}p\big],\\ &\bar{y}(g(\bar{x},\bar{y})-\theta_0)p-\frac{\sigma^2}{2}\frac{\partial}{\partial \bar{y}}\big[\bar{y}^2p\big]-\sigma^2\frac{\partial}{\partial \bar{x}}\big[\bar{x}\bar{y}p\big]\bigg),
	\end{split}
	\end{align}
	
	\noindent The no flux or reflecting boundary conditions are then
	
	\begin{equation}
	\bf{n}\cdot\bf{J}=0.
	\end{equation}
	
	\noindent where $\bf{n}$ is the outward normal to the surface. In our case this condition is automatically satisfied at the $\bar{x}=0,\bar{y}=0$ boundaries and we only need to impose it on the diagonal line as well as the $\bar{x}=3,\bar{y}=3$ boundaries. 
	
	\subsection{Adding death rates to the model}
	
	To explore what happens in the case of adding death rates to the model we first need to take a step back and present the dimensional system.
	The growth rates for two microbial population are the combined growth rates, given by the Monod model of uninhibited growth, minus the death rate ($d_i$):
	
	\begin{align}
	\mu_1(s)=\frac{\mu_m^1s}{K_s^1+s}-d_1\\
	\mu_2(s)=\frac{\mu_m^2s}{K_s^2+s}-d_2
	\end{align}
	
	\noindent Here, $\mu_m^i$ is the maximum specific growth rate, $K_s^i$ is the Michaelis constant in the specific growth rate expression and $s$ is the substrate concentration.\\
	
	We are interested in the points of intersection between these two curves, which are always two since the equation $\mu_1(s)=\mu_2(s)$ is quadratic in $s$. In the absence of death rate it is easy to deduce that there can potentially be two intersection points with one being always zero and the other for positive $s$ and positive values of the growth rates. Instead of the positive intersection point it is possible that we have one for negative $s$ value and negative growth rates but this crossing is no longer of any biological interest. So we have only two interesting cases, (a) one intersection point at zero or (b) two points, one is zero and the other somewhere in the upper right quadrant. For (b) it is trivial to find the relations between the parameters of the two growth rates depending on the assumption you make about the growth rates, i.e. which is largest as $s \rightarrow \infty$.

	\begin{figure}[H]
		\begin{center}
			\subfigure[Growth rates without death rate.]{
				\resizebox*{6cm}{!}{\includegraphics{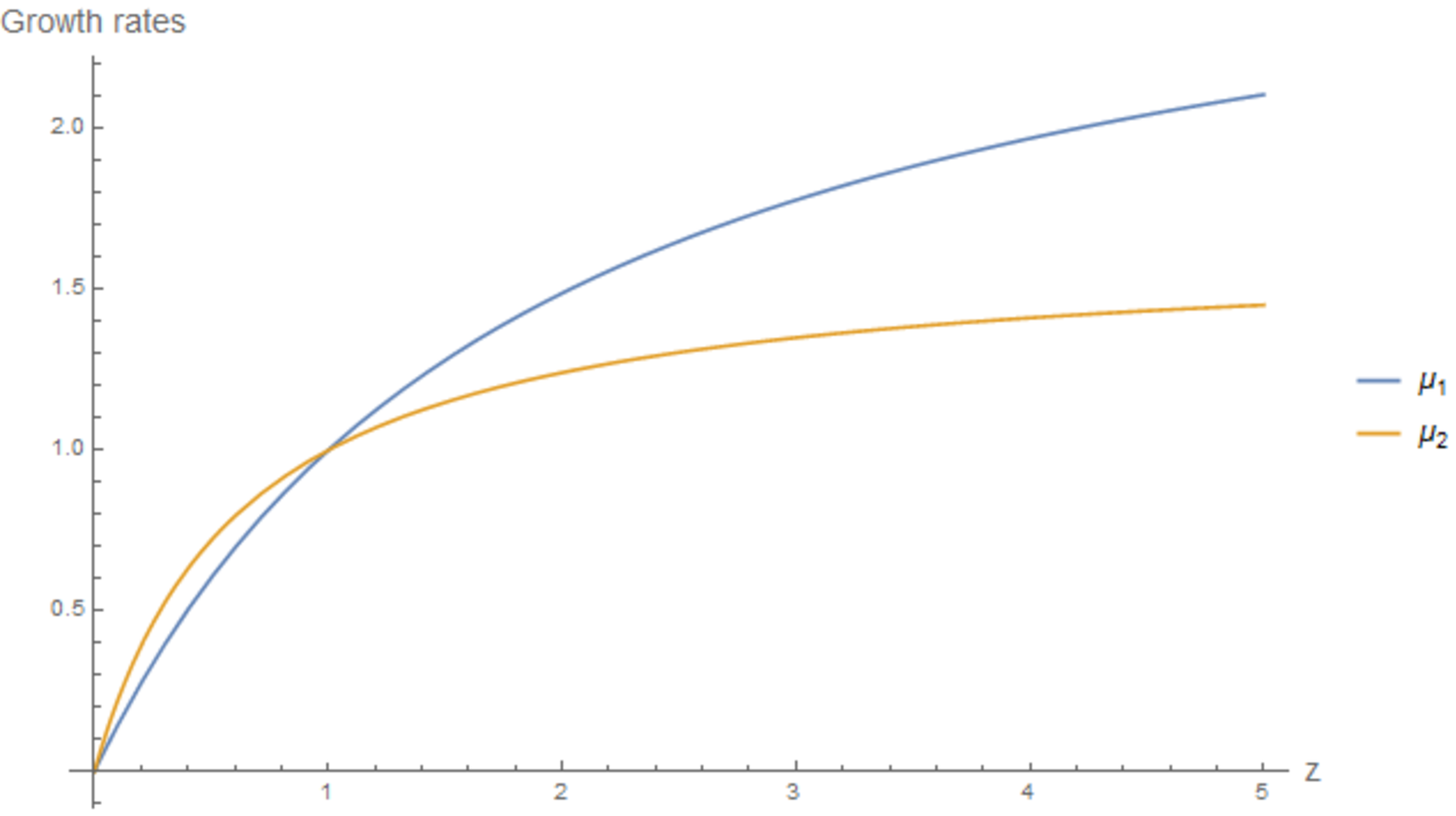}}}\hspace{5pt}
			\subfigure[Growth rates with death rate.]{
				\resizebox*{6cm}{!}{\includegraphics{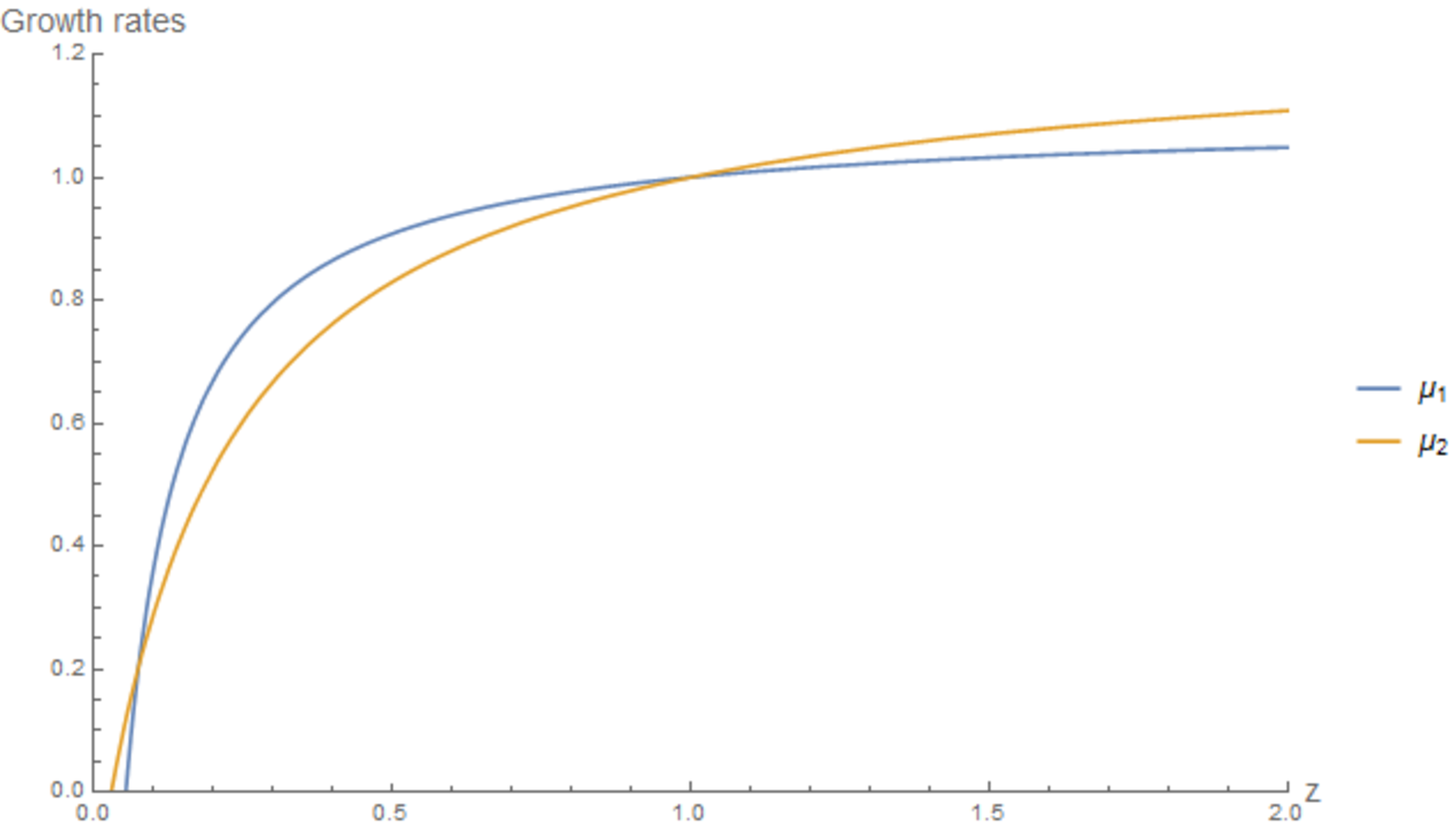}}}
			\caption{Plot of the dimensional growth rates with and without death rates for different substrate values.}
		\end{center}
	\end{figure}
	
	On the other hand, including two distinct death rates complicates the situation. Now there are three biologically relevant cases which give very different relations between the parameters, (a) both the intersection points being at the upper right quadrant, (b) the first being at the lower right and the second on the upper right and (c) the first being in the lower left quadrant and the second in the upper right. In all the other cases both intersections points occur at negative growth rate values.\\
	
	First we will focus on (b) and (c) since they are the simplest cases and admit only one interesting point of intersection. In the analysis that follows it is always assumed that $\mu_1(s)>\mu_2(s)$ as $s \rightarrow \infty$, and hence obtaining the first relation:
	
	\begin{equation}
	\mu_m^1-d_1 > \mu_m^2-d_2
	\end{equation}
	
	Moreover, we want the end-values of the growth rates to be positive which gives the next two relations:
	
	\begin{align}
	\mu_m^1>d_1\\
	\mu_m^2>d_2
	\end{align}
	
	The simplest case is (b) where both intersection points are for $s>0$ and the first admits a negative growth rate whereas the second a positive. Let us have a look at the additional conditions imposed by (b). The solutions to the quadratic equation $f(s)=g(s)$ are:
	
	\begin{equation*}
	s_{1,2}=\frac{-B \mp \sqrt{B^2-4AC}}{2A}
	\end{equation*}
	
	\noindent Where A, B and C are as follows:
	
	\begin{align*}
	A&=\mu_1-d_1+d_2-\mu_2\\
	B&=K_2(\mu_1-d_1+d_2)-K_1(\mu_2-d_2+d_1)\\
	C&=K_1K_2(d_2-d_1)
	\end{align*}
	
	In order for both $s_1, s_2$ to be positive we just need the smaller, $s_1$, to be positive. For that to happen we need the following conditions, (I) $B<0$, (II) $C>0$ so that $\sqrt{B^2-4AC}<-B$, since we already know from (48) that $A>0$, and finally (III) $B^2>4AC$ in order to have real solutions. From these conditions only one is of immediate importance and that is (II) due to the fact that (I) and (III) will be satisfied by another stronger condition we need in order to have the desirable geometry ($\mu_1(s_2)>0, \mu_1(s_1)<0$). So, $C>0$ implies:
	
	\begin{equation}
	d_2>d_1
	\end{equation}
	
	In addition if we also require $\mu_1(s_1)=\mu_2(s_1)<0$ and $\mu_1(s_2)=\mu_2(s_2)>0$ we get one final condition:
	
	\begin{equation}
	\frac{K_1d_1}{\mu_1-d_1}>\frac{K_2d_2}{\mu_2-d_2}
	\end{equation}
	
	What that inequality says is that $\mu_2$ crosses the s-axes before $\mu_1$, meaning that $s_0^1>s_0^2$ where, $\mu_1(s_0^1)=\mu_2(s_0^2)=0$. Using (48) and (51) another condition is obtained which is weaker than (52), i.e. necessary but not sufficient, but useful for the stability analysis that follows.
	
	\begin{equation}
	K_1>K_2
	\end{equation}
	
	Inequality (52) was found using Mathematica and we can again using Mathematica show that if it is satisfied then $B<0$ and $B^2-4AC \geq 0$. Hence, the only two condition needed for case (b) are (51) and (52).\\
	
	The nondimensionalisation as well as the results of the stability analysis are found in Appendix B. The linear stability is summed up in the table below:
	\bigskip
	
	\begin{table}
		\caption{Stability analysis w/ death rate.}
		\centering
		{\begin{tabular}[l]{@{}lcc}\toprule                    
				\textbf{Steady-State} & \textbf{Conditions for $\lambda_1<0$} & \textbf{Conditions for $\lambda_2<0$}\\ \hline
				$x=y=0$ & $f(z_f)>\theta$ & $g(z_f)>\theta$ \\
				$x=0, g(z)=\theta$ & $g(z_f)>\theta$  & $\theta<1$  \\
				$y=0, f(z)=\theta$ & $f(z_f)>\theta$  & $\theta>1$ \\ 
				$f(z)=g(z)=\theta$ & $z_f>1$ AND $\theta=1$ &NA $\lambda_2=0$ \\
				\hline   
		\end{tabular}}
	\end{table}
	
	\bigskip
	
	If now case (c) is considered, namely the case where one intersection point is at the upper right quadrant and the other in the lower left, the exact same stability analysis results are recovered as was expected. In case (c) there is one condition opposite to (51) which, when combined with (48) admits one extra condition:
	
	\begin{align}
	d_1>d_2\\
	a_1>a_2
	\end{align}
	
	\noindent Other than that the stability analysis remains the same as in Table\\
	
	As stated before this analysis is valid for cases (b) and (c) when there is a single intersection point in the upper right quadrant. In case (a)  where there are two nondimensionalisation and linear stability analysis is much more complicated and very difficult to perform. So instead we summarize the stability of the system via numerical simulation in Figure 9. 
	
	\begin{figure}[H]
		\begin{center}
			\resizebox*{12cm}{!}{\includegraphics{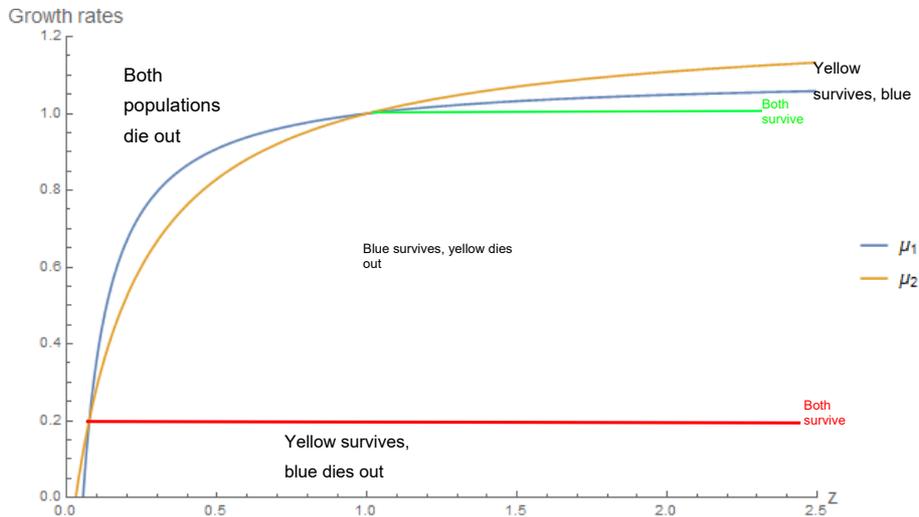}}
			\caption{Stability of system with death rates and two intersection points in the upper quadrant with initial populations being the same.} \label{fig1}
		\end{center}
	\end{figure}
	
	We can deduce that the stability seems to be exactly what we expected according to the cases where we had one intersection point, i.e. the population surviving in each segment of the parameter space is the one with the highest growth rate.

	\section{Results}
	
	Due to the fact that $z_f$ is usually large we chose an arbitrarily large value equal to 15000. Results for the full system are not affected by the value of $z_f$, except as regards numerical accuracy,while our asymptotic analysis and further simulation of the simplified system and Fokker-Planck take $z_f$ asymptotically large.
	
	\subsection{Simulation of Langevin equations for full system without death}
	
	To simulate the dynamics of the full system close to the line of steady states we use the Euler-Maruyama method on the system of SDEs (6-8) and (9-12). The initial conditions for $x,y$ are the given by $y=z_f-x-1$ which is the line of steady states and specifically we chose $x=z_f/2 -1,y=z_f/2$. Moreover $z$ always starts from the value of 1 so that the system is initiated on the line of steady states. \\
	
	The Euler-Maruyama method has a strong convergence of order 1/2. Alternatively we can use the Milstein method with a strong convergence of order 1 \cite{higham_algorithmic_2001}. We noticed no difference in the numerical results from the two approaches and hence we performed the simulations using the former. In the following figures we present the evolution of the two populations for different values of $\theta$ and different noise intensities.\\
	
	\subsubsection{Evolution of  microbial populations for general noise}
	For the simulation of the general noise case we initially picked the values of the noise intensities such that they are at the same order as the the ones used at the dilution rate induced case. Meaning that for $\sigma_3$ we picked high intensity to account for the fact that environmental stochasticity affects z dramatically because of the $z_f$ term in the equation. In addition we increased the values of $\sigma_1,\sigma_2$ to see if there is any qualitative difference to the dynamics. We can do that since the noise terms are different for $x,y$ and $z$. For $\theta=0.99$ though we can see that when the noise is high the competition persists for long and although $y$ is the only one surviving there could be values of the intensity that either does. We do not observe something similar for $\theta=1.01$ where it is again clear who the winner is and the results agree with the deterministic case. So we further see an advantage of the $x$ population of the $y$.
	
	\begin{figure}[H]
		\begin{center}
			\subfigure[$\theta=1,\sigma_1=0.0006,\sigma_2=0.0007,\sigma_3=9$.]{
				\resizebox*{6cm}{!}{\includegraphics{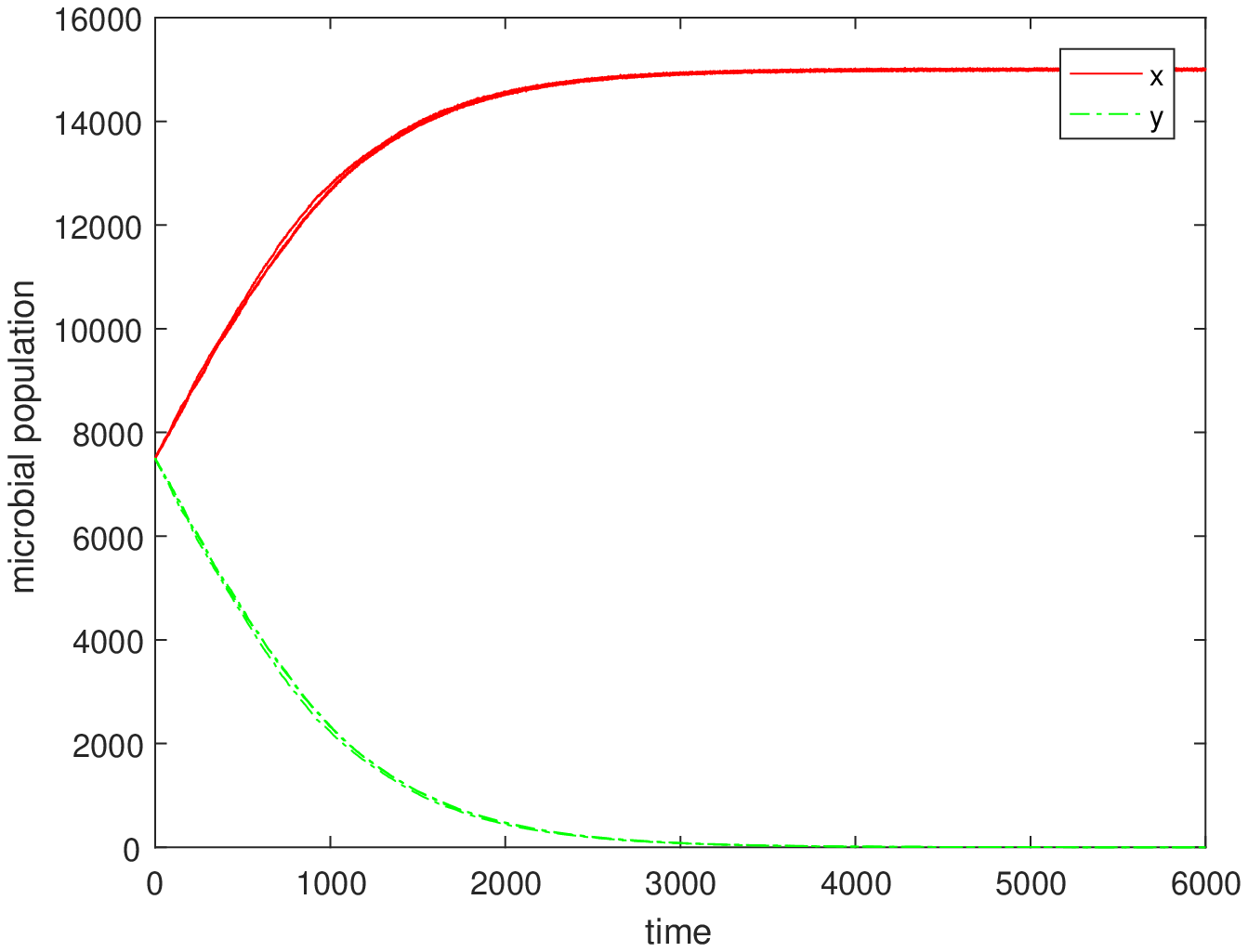}}}\hspace{5pt}
			\subfigure[$\theta=1,\sigma_1=0.05,\sigma_2=0.07,\sigma_3=9$.]{
				\resizebox*{6cm}{!}{\includegraphics{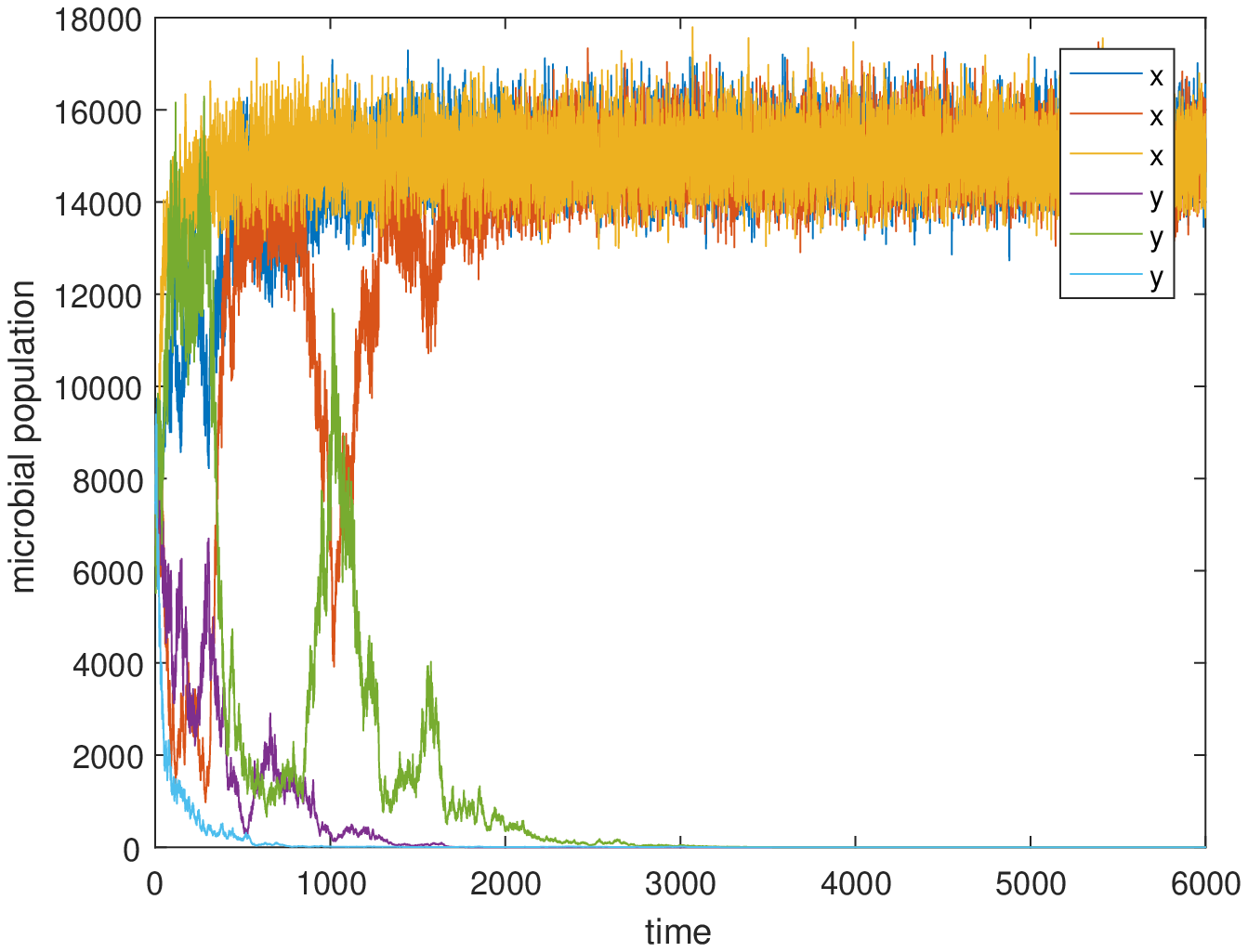}}}
			
			\subfigure[$\theta=1.02,\sigma_1=0.0006,\sigma_2=0.0007,\sigma_3=9$.]{
				\resizebox*{6cm}{!}{\includegraphics{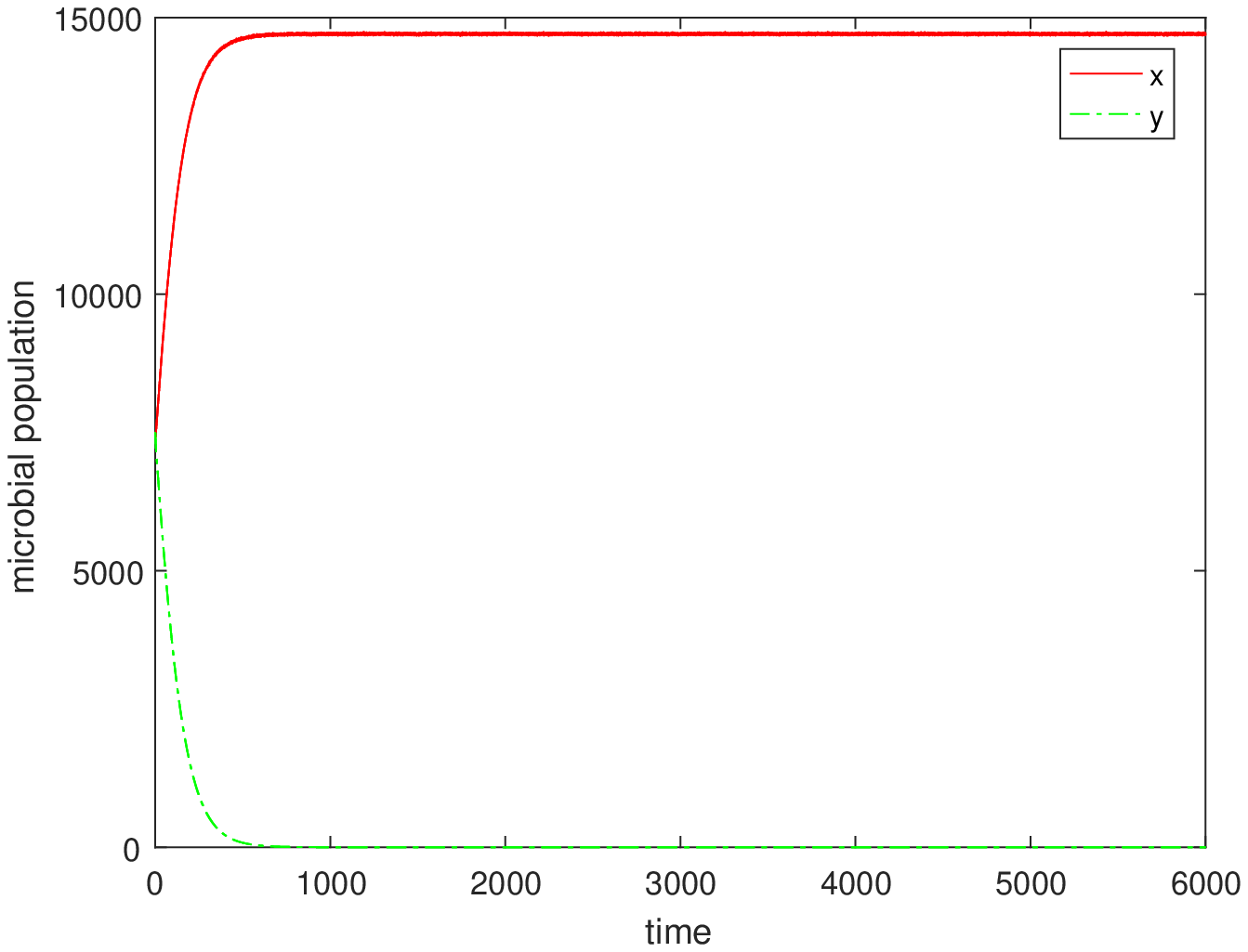}}}\hspace{5pt}
			\subfigure[$\theta=1.02,\sigma_1=0.05,\sigma_2=0.07,\sigma_3=9$.]{
				\resizebox*{6cm}{!}{\includegraphics{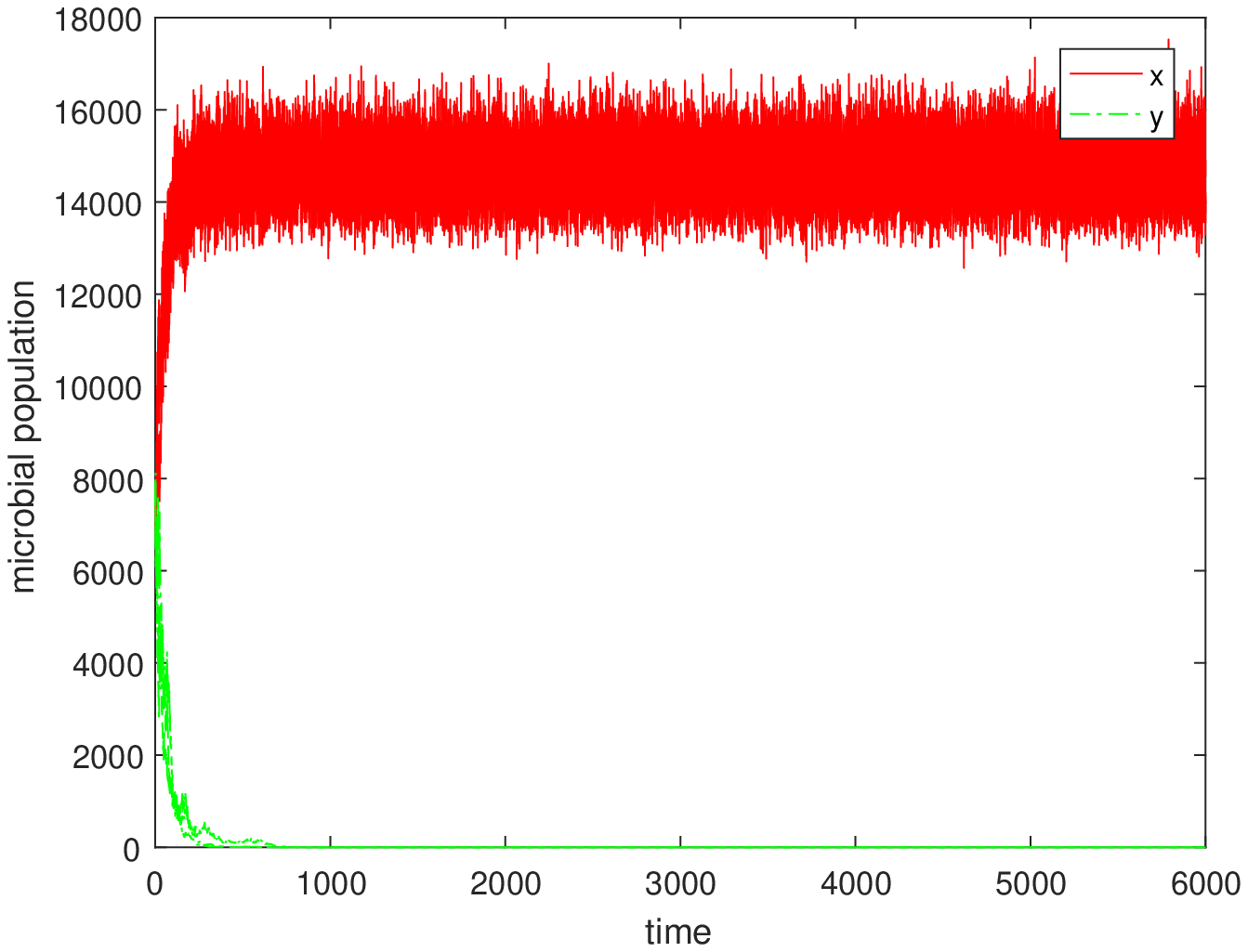}}}

			\subfigure[$\theta=0.98,\sigma_1=0.0006,\sigma_2=0.0007,\sigma_3=9$.]{
				\resizebox*{6cm}{!}{\includegraphics{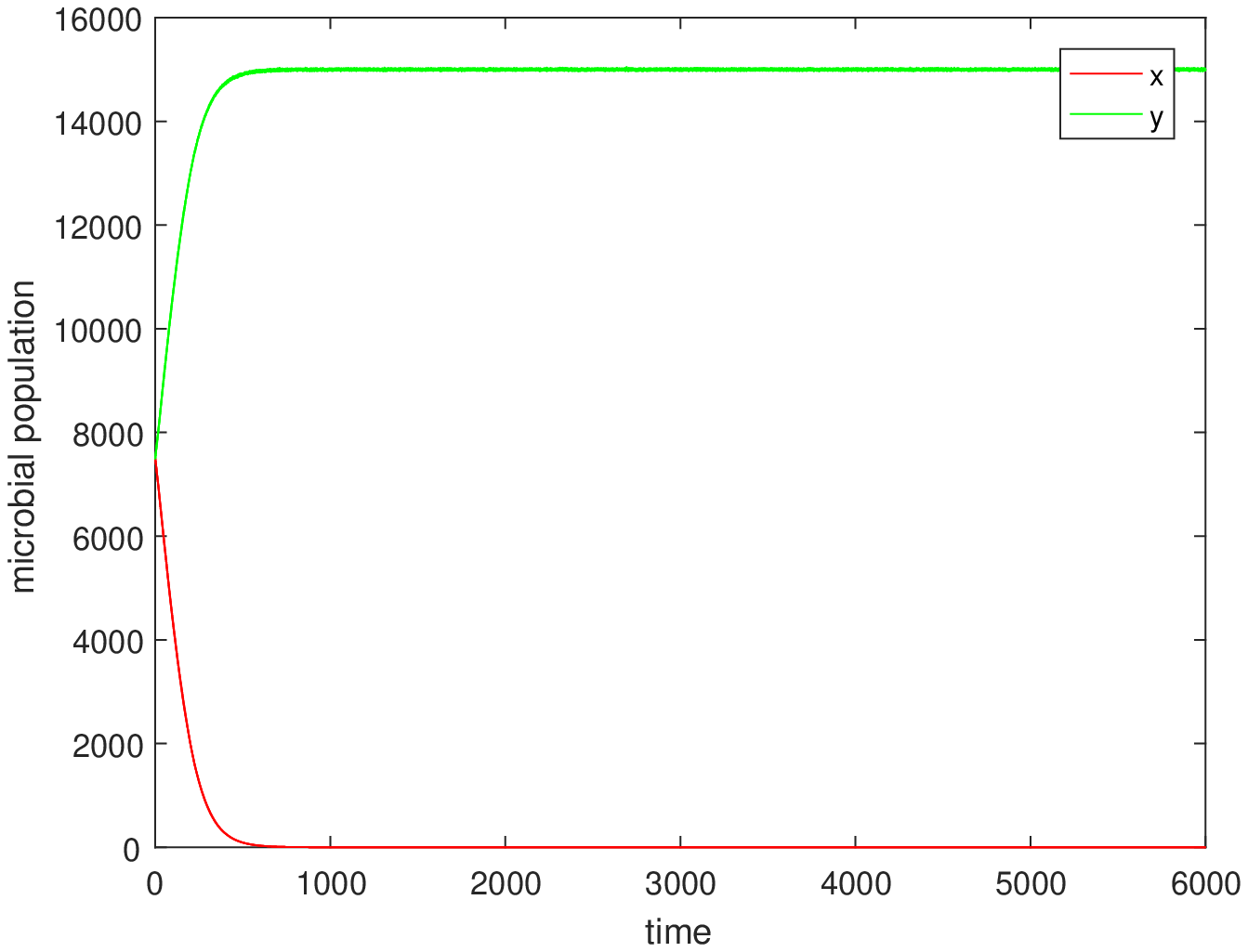}}}\hspace{5pt}
			\subfigure[$\theta=0.98,\sigma_1=0.05,\sigma_2=0.07,\sigma_3=9$.]{
				\resizebox*{6cm}{!}{\includegraphics{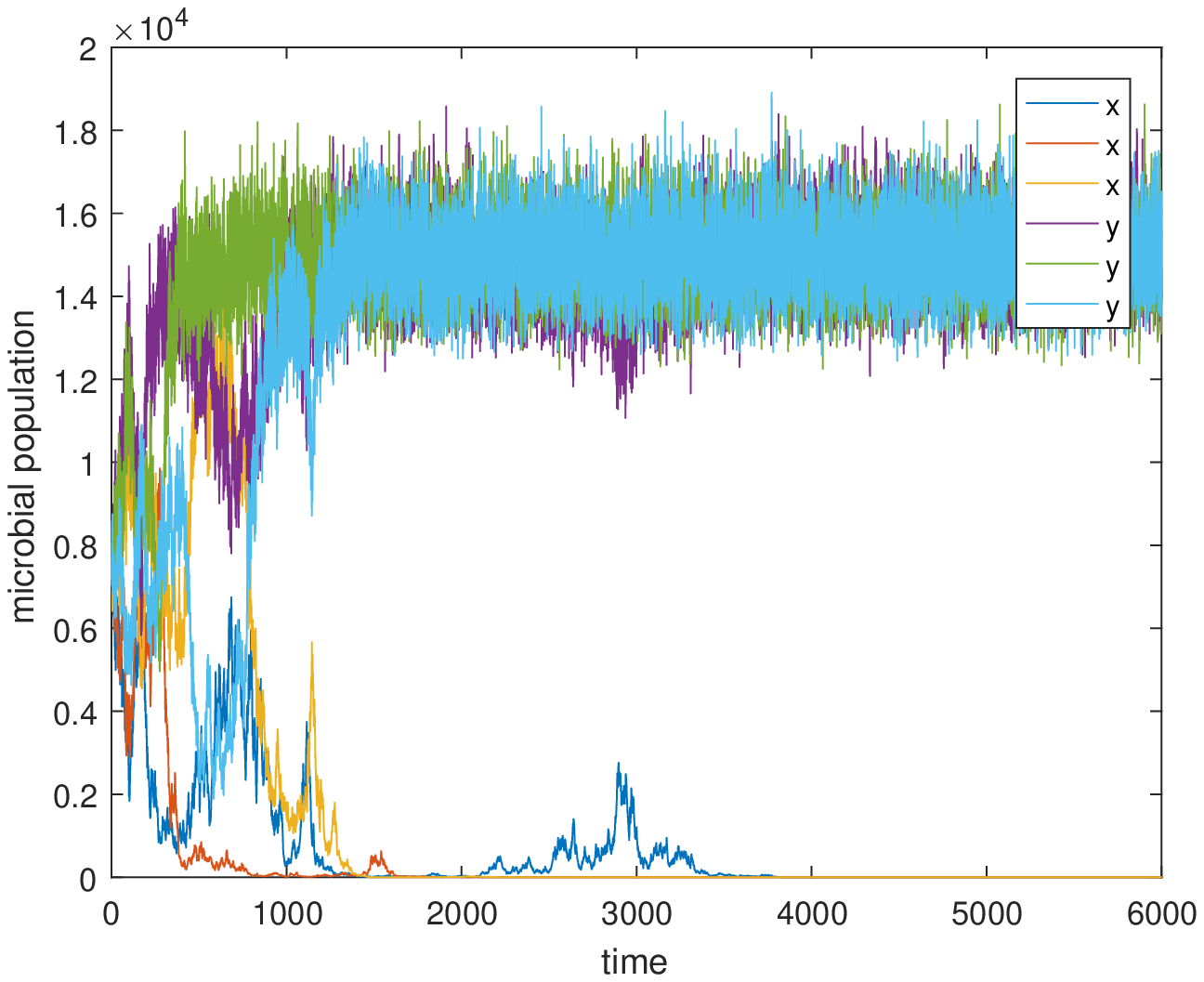}}}
			\caption{Plot of 3 runs of non-dimensional microbial populations versus non-dimensional time for varying values of the dilution rate and noise intensity.}
		\end{center}
	\end{figure}
	
	\begin{figure}[H]
		\begin{center}
			\subfigure[$\theta=0.99,\sigma_1=0.0006,\sigma_2=0.0007,\sigma_3=9$.]{
				\resizebox*{6cm}{!}{\includegraphics{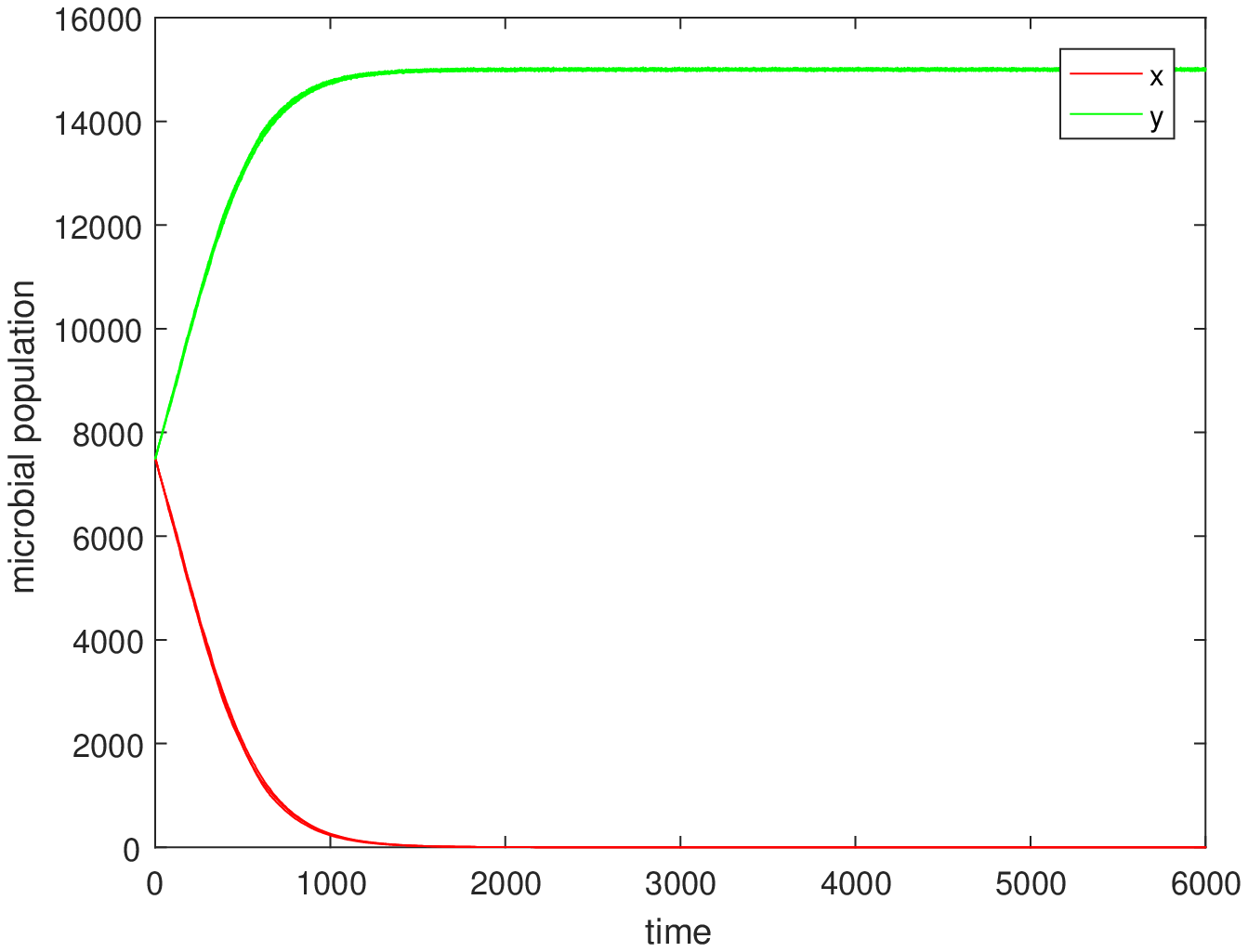}}}\hspace{5pt}
			\subfigure[$\theta=0.99,\sigma_1=0.05,\sigma_2=0.07,\sigma_3=9$.]{
				\resizebox*{6cm}{!}{\includegraphics{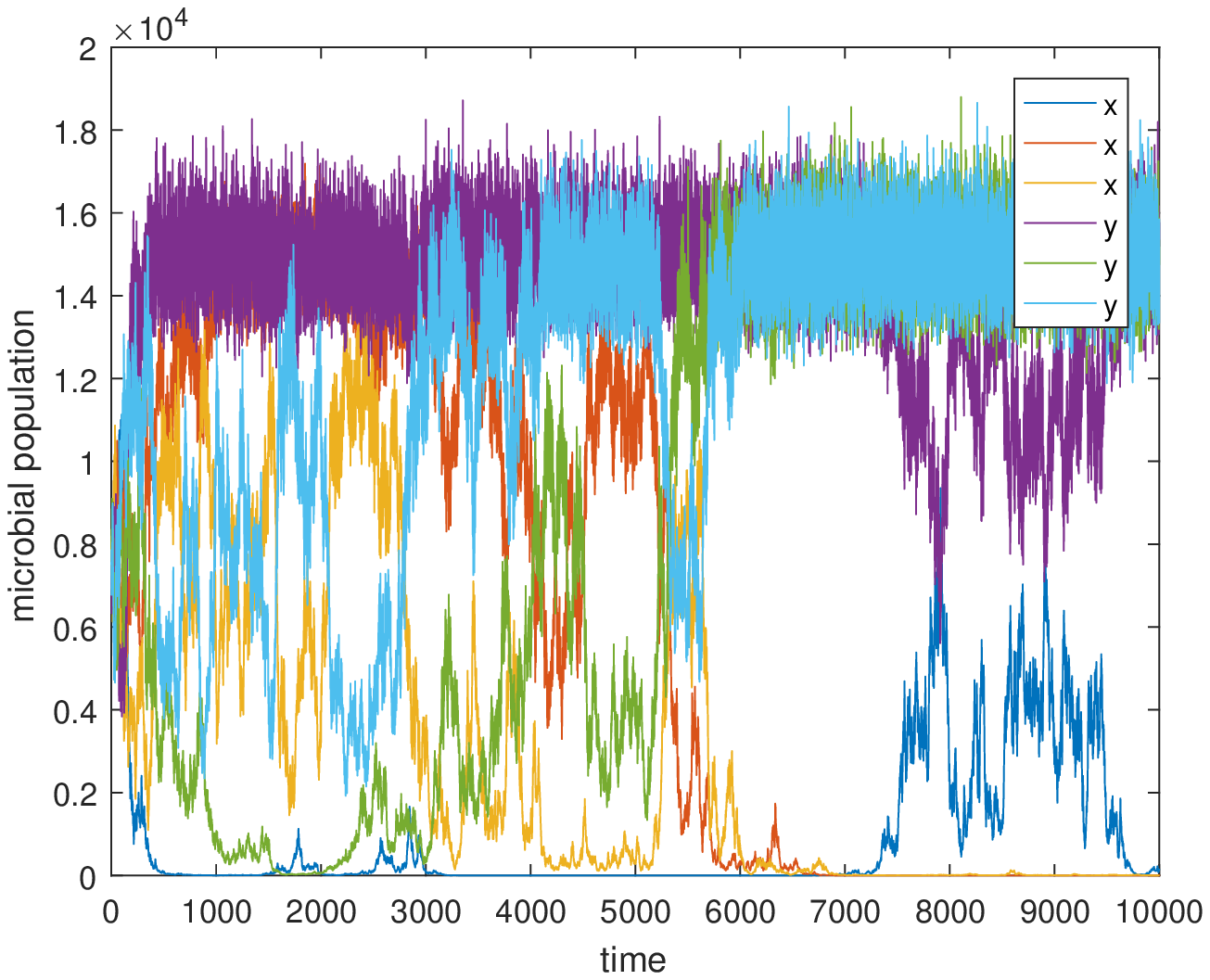}}}
			
			\subfigure[$\theta=1.01,\sigma_1=0.0006,\sigma_2=0.0007,\sigma_3=9$.]{
				\resizebox*{6cm}{!}{\includegraphics{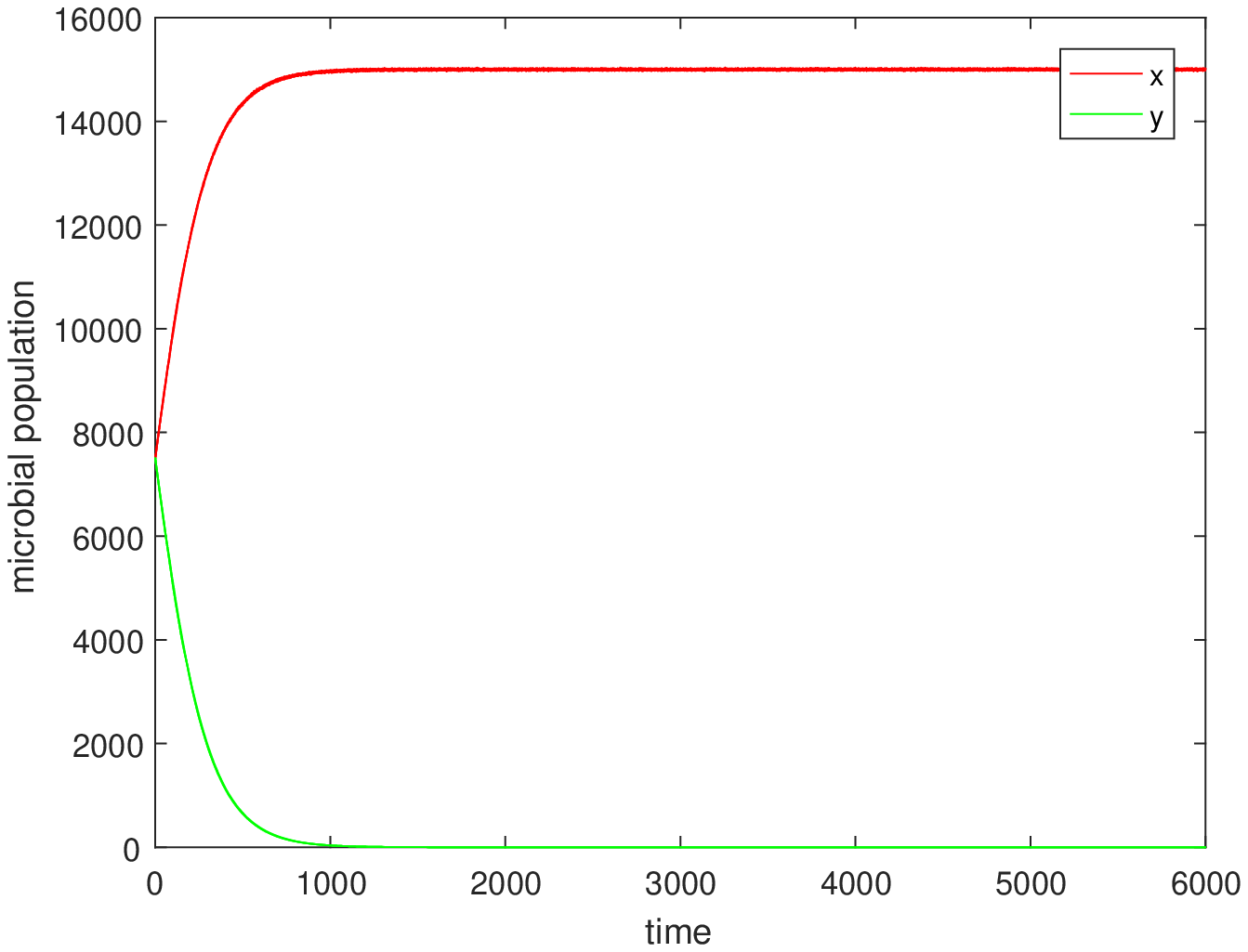}}}\hspace{5pt}
			\subfigure[$\theta=1.01,\sigma_1=0.05,\sigma_2=0.07,\sigma_3=9$.]{
				\resizebox*{6cm}{!}{\includegraphics{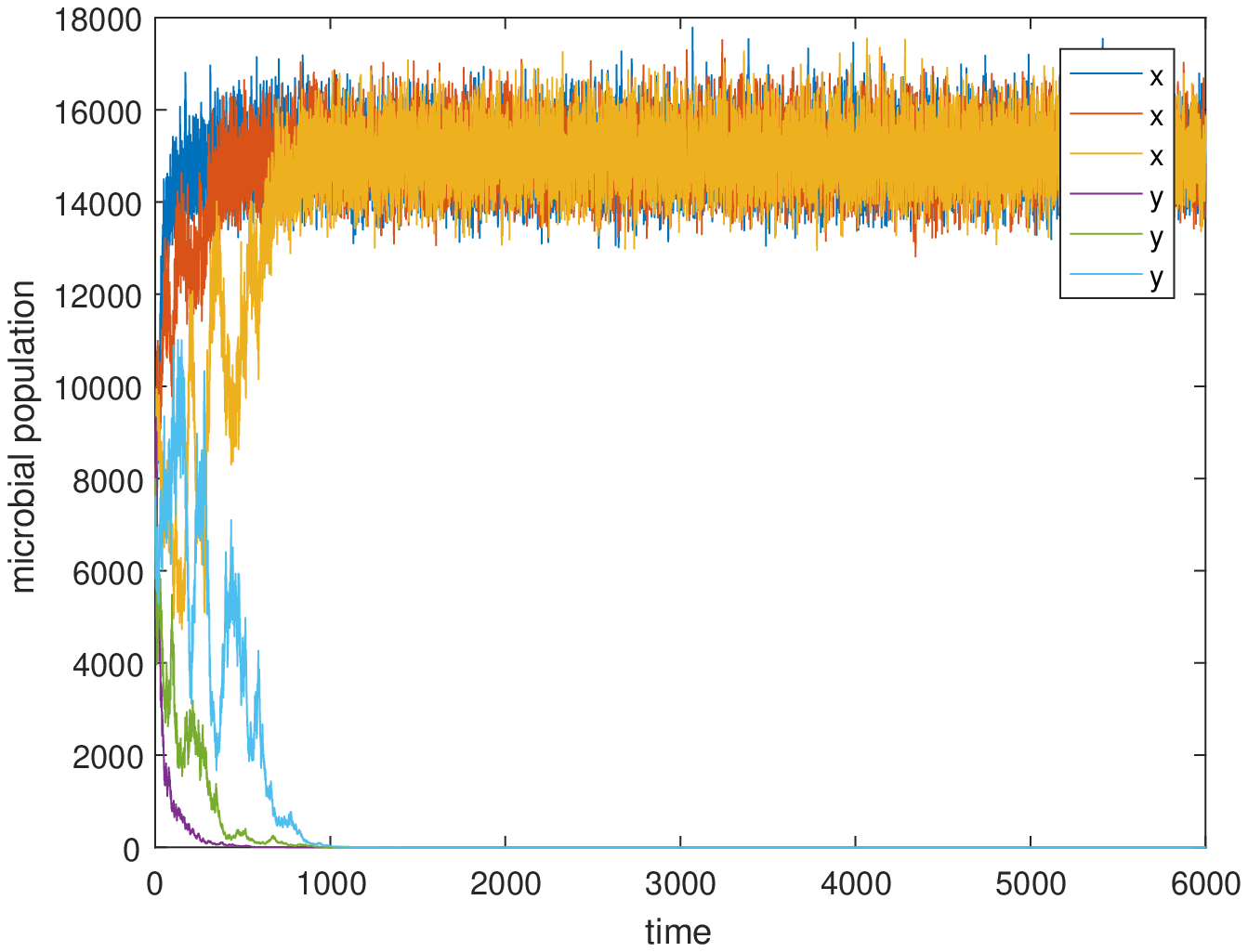}}}
			\caption{Plot of 3 runs of non-dimensional microbial populations versus non-dimensional time for varying values of the dilution rate and noise intensity.}
		\end{center}
	\end{figure}

	It is interesting to explore what happens to the case of high noise in $x,y$ and decreased noise in $z$
	
	\begin{figure}[H]
		\begin{center}
			\subfigure[$\theta=1,\sigma_1=0.05,\sigma_2=0.07,\sigma_3=0.01$.]{
				\resizebox*{6cm}{!}{\includegraphics{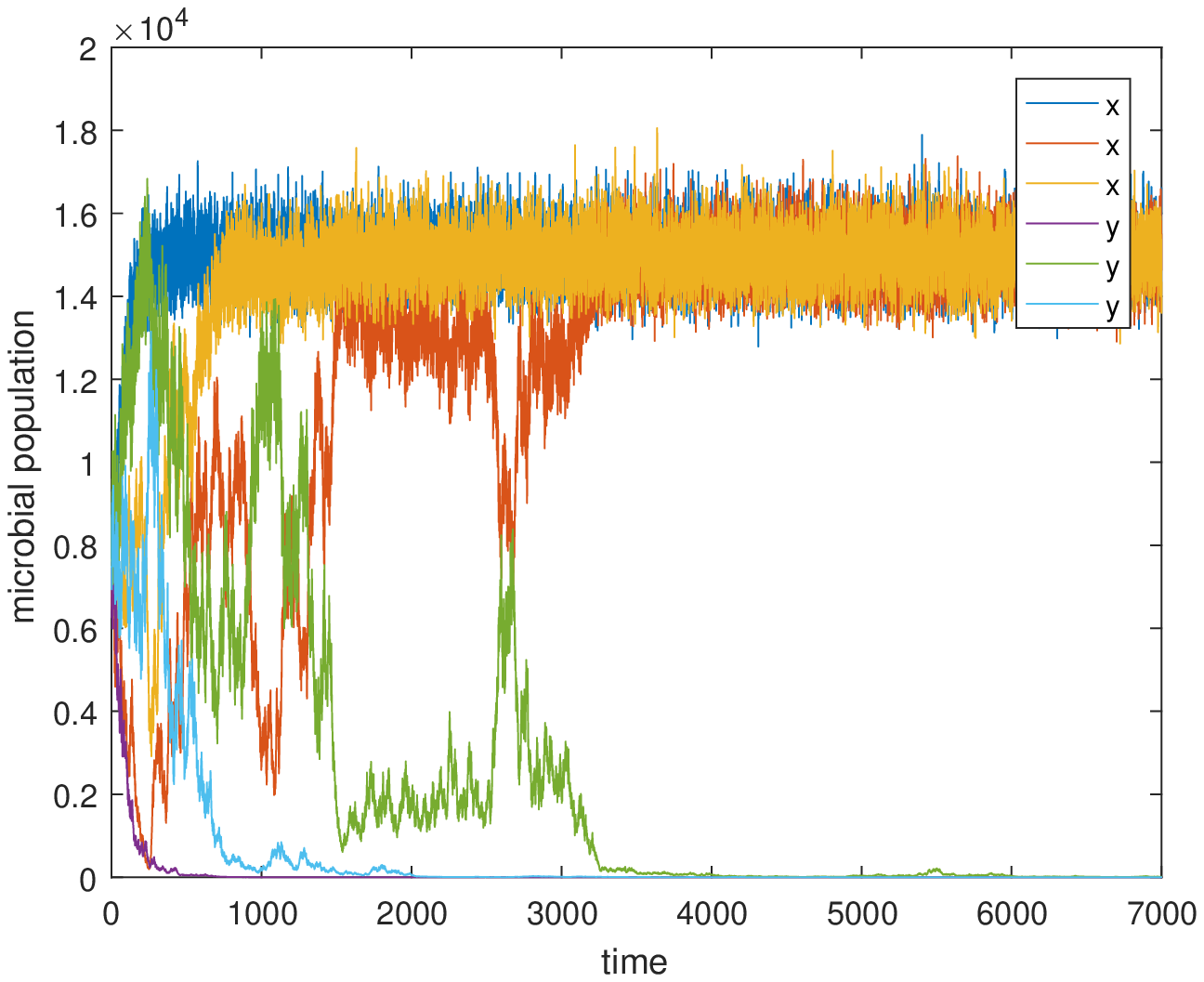}}}\hspace{5pt}
			\subfigure[$\theta=1,\sigma_1=0.05,\sigma_2=0.07,\sigma_3=0.00001$.]{
				\resizebox*{6cm}{!}{\includegraphics{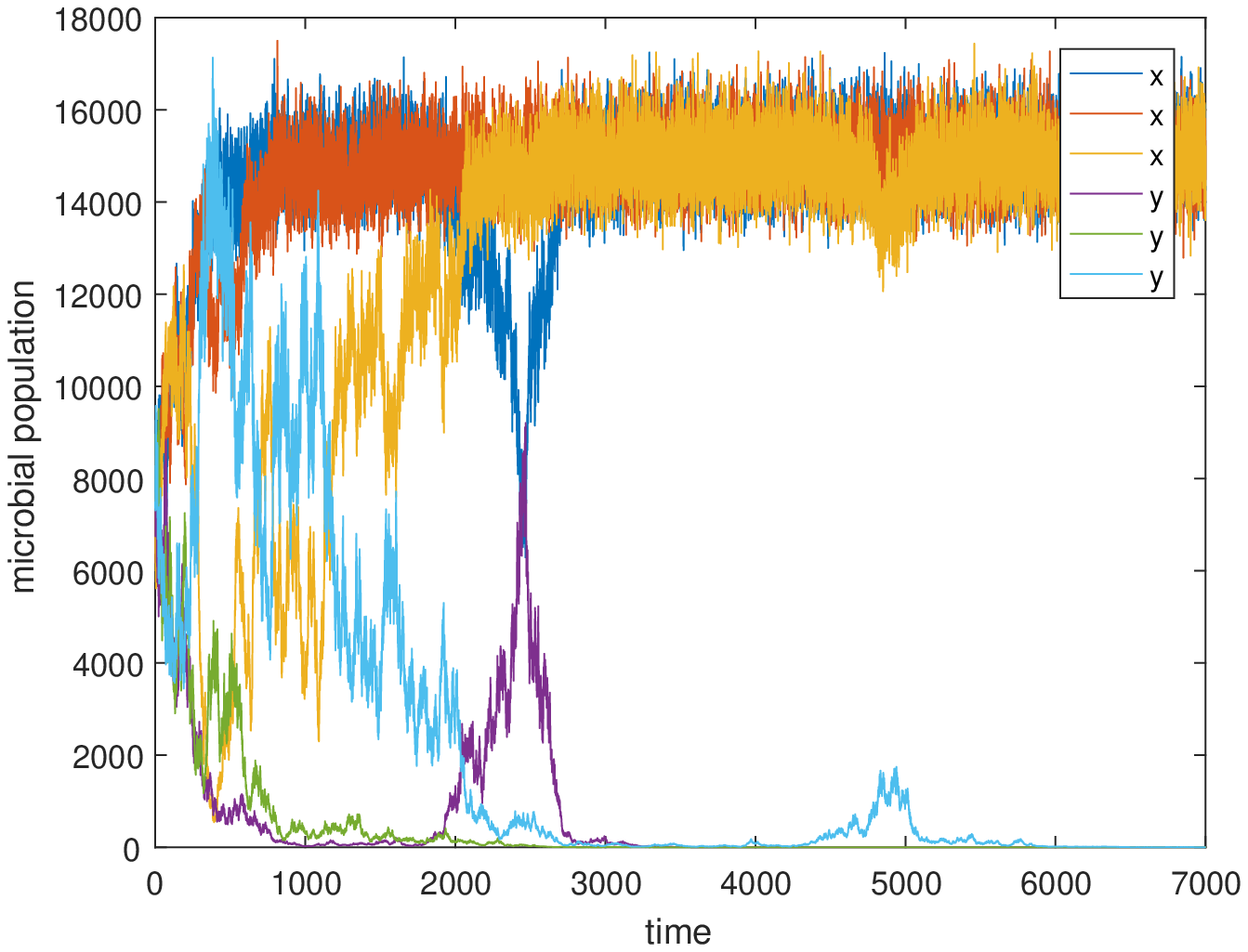}}}
			\caption{Plot of 3 runs of non-dimensional microbial populations versus non-dimensional time for varying values the $\sigma_3$ noise intensity.}
		\end{center}
	\end{figure}
	
	The two plots combined with Figure 9(b) show that the variation the dynamics of $x,Y$ are mainly affected by the stochasticity in $x,y$ and that noise in $z$ has little to no effect whether $\sigma_3$ is larger or smaller than the noise intensity. That can justify the replacement of the differential equation for $z$ with the algebraic equation found using the asymptotic analysis when the system is close to that line.\\
	
	\subsubsection{Evolution of  microbial populations for dilution rate induced noise}
	
	 The plots of Figure 12 and 13 show a clear dominance of the x population for the case where the mean value of the dilution rate is 1. We can further see that higher noise favours the population with the highest maximum growth rate whereas lower favours the other. This is particularity evident in the last two plots. For $\theta=0.99$ deterministically we are in the region were y should survive and x must vanish and with low noise we see that this is the case. However as noise intensifies we can see that there is a longer competition between the two populations. For particular values of $\theta,\sigma$ we can create a situation where both populations coexist for a very long period were a very extended simulation needs to take place to determine which one will eventually survive. It is also worth noticing that there is a contrasting effect of the noise intensity depending on which region of $\theta$ values the system is in. If it is $\geq 1$ then higher intensity drives the system to a steady state faster whereas if $\leq 1$ higher intensity has the opposite effect.\\
	
	The plots for the case of $\theta=1$ are in contrast to what was found in \cite{stephanopoulos_stochastic_1979} where because of the fact that z was considered constant, the drift term disappeared in the Fokker-Planck equation and the evolution was purely stochastic. As a result both populations had similar chances of survival if they had the same initial conditions. This is not the case here.
	
	\begin{figure}[H]
		\begin{center}
			\subfigure[$\theta=1,\sigma=0.0006$.]{
				\resizebox*{6cm}{!}{\includegraphics{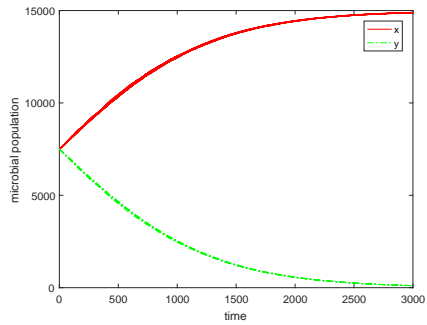}}}\hspace{5pt}
			\subfigure[$\theta=1,\sigma=0.001$.]{
				\resizebox*{6cm}{!}{\includegraphics{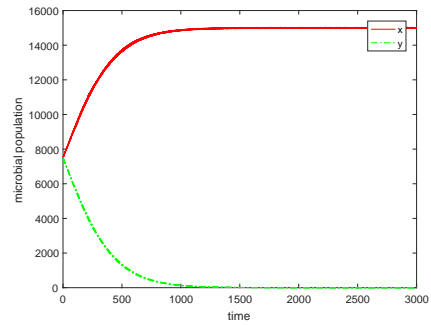}}}
			
			\subfigure[$\theta=1.02,\sigma=0.0006$.]{
				\resizebox*{6cm}{!}{\includegraphics{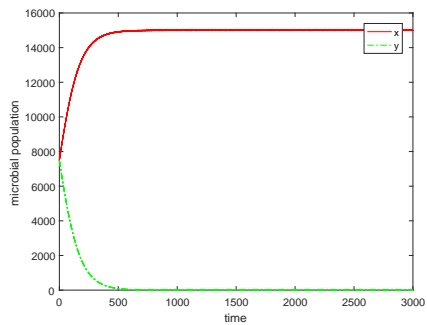}}}\hspace{5pt}
			\subfigure[$\theta=1.02,\sigma=0.001$.]{
				\resizebox*{6cm}{!}{\includegraphics{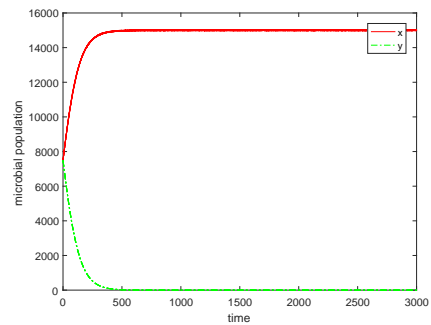}}}
			\caption{Plot of 10 runs of non-dimensional microbial populations versus non-dimensional time for varying values of the dilution rate and noise intensity.}
		\end{center}
	\end{figure}
	
	\begin{figure}[H]
		\begin{center}
			\subfigure[$\theta=0.98,\sigma=0.0006$.]{
				\resizebox*{6cm}{!}{\includegraphics{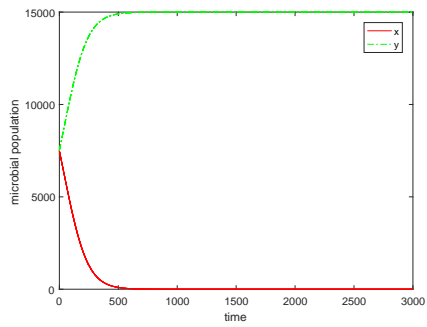}}}\hspace{5pt}
			\subfigure[$\theta=0.98,\sigma=0.001$.]{
				\resizebox*{6cm}{!}{\includegraphics{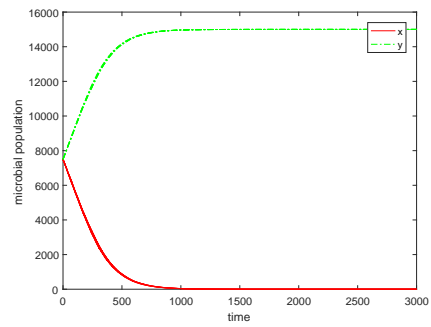}}}
			
			\subfigure[$\theta=0.99,\sigma=0.0006$.]{
				\resizebox*{6cm}{!}{\includegraphics{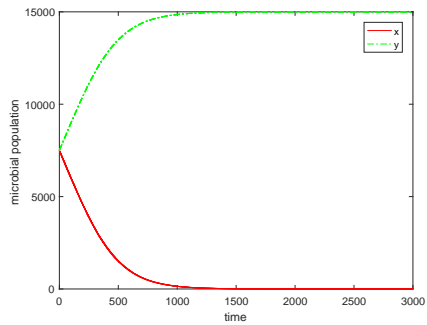}}}\hspace{5pt}
			\subfigure[$\theta=0.99,\sigma=0.001$.]{
				\resizebox*{6cm}{!}{\includegraphics{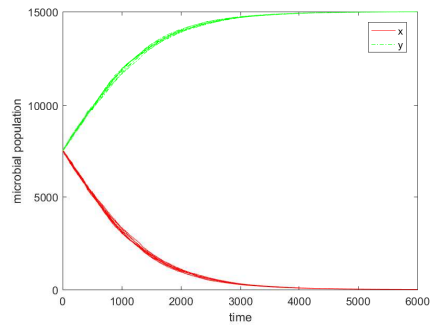}}}
			\caption{Plot of 10 runs of non-dimensional microbial populations versus non-dimensional time for varying values of the dilution rate and noise intensity.}
		\end{center}
	\end{figure}
	
	Despite the fact that the noise term is the same in all three equations of the full system it is interesting to see what happens if we kept the noise intensity, $\sigma$, the same in (10) and (11) but decrease it in (12).
	
	\begin{figure}[H]
		\begin{center}
			\subfigure[$\theta=1,\sigma_{x,y}=0.0006, \sigma_{z}=0.000006$.]{
				\resizebox*{10cm}{!}{\includegraphics{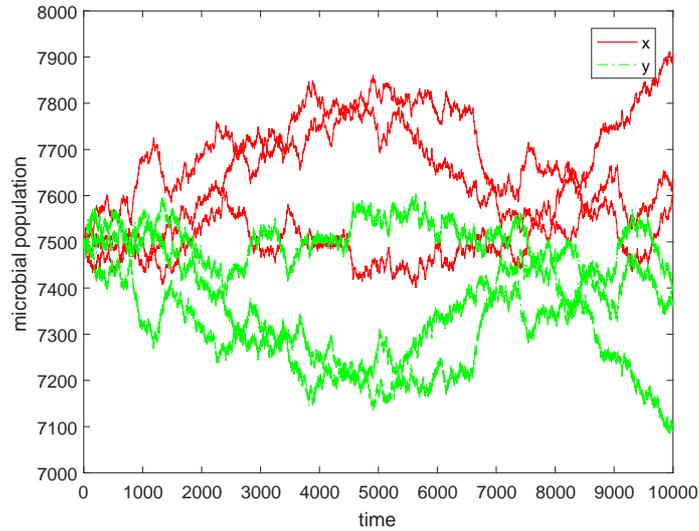}}}\hspace{5pt}
			\subfigure[$\theta=0.98,\sigma_{x,y}=0.03, \sigma_{z}=0.000006$.]{
				\resizebox*{10cm}{!}{\includegraphics{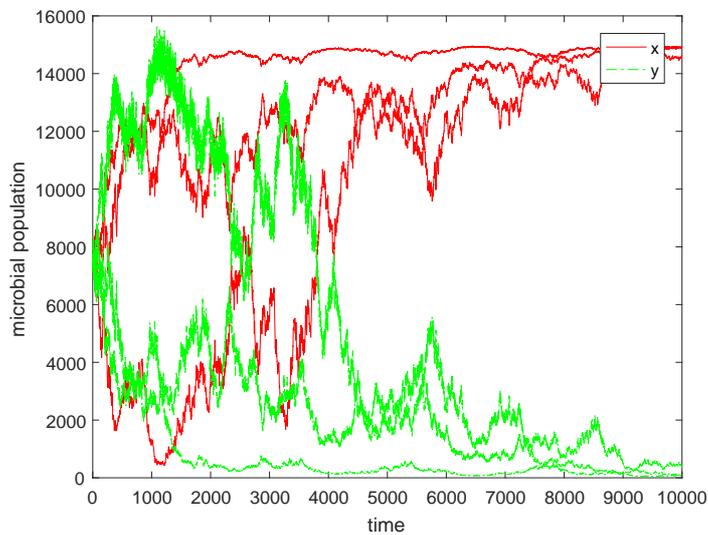}}}
			\caption{Plot of 3 runs of non-dimensional microbial populations versus non-dimensional time for different values of the noise intensity in the x,t stochastic equations and the z equation.}
		\end{center}
	\end{figure}
	
	We can clearly see that the lower noise in $z$ causes the dynamics to become more noisy and hence reach a steady-state value at significantly larger times than the case were the noise of z is greater.

	\subsection{Simulation of Differential-Algebraic system for general noise}
	
	Having justified the use of the reduced system for the case of general noise we can use the derived Fokker-Planck equation (39) to simulate the probability density function of the system. The advantage of using the DAE system without noise in the algebraic equation is that if we solve for z as mentioned in the previous section and substitute into the SDEs for $x,y$ we can then derive a Fokker-Planck equation in two "spatial" dimensions and time which allows us to explore the evolution of the system for low noise and long times with much lower computational cost, higher accuracy and more clarity. We may also refer to the spatially 3D Fokker-Planck system which has a single derivative in t and double derivatives in x, y, z. This system also involves a single z-derivative of the algebraic term in square brackets on the right-hand side of equation (8). Our emphasis guided by the numerical studies is on the solution of the system when the latter (algebraic) term is zero. Setting that term to zero then creates a specific connection between x, y, z derivatives. This allows any z derivative for example to be replaced by a combination of x, y derivatives and that in turn leads to the Fokker-Planck system becoming one in t, x, y only as in the present working. The reason for not wanting to use a noise intensity that is very low is the fact that the smaller the noise intensity the smaller the diffusion term of the Fokker-Planck equation and the finer the mesh needs to be to solve the equations correctly. This 2D Fokker-Planck was numerically solved using the finite element method (FEM) in Mathematica 11.\\
	
	\begin{figure}[H]
		\begin{center}
			\subfigure[FP equation for $t=1$ (Green), $t=500$ (Red), $t=10000$ (Yellow), $\sigma_1=0.05$ and $\sigma_2=0.07$.]{
				\resizebox*{10cm}{!}{\includegraphics{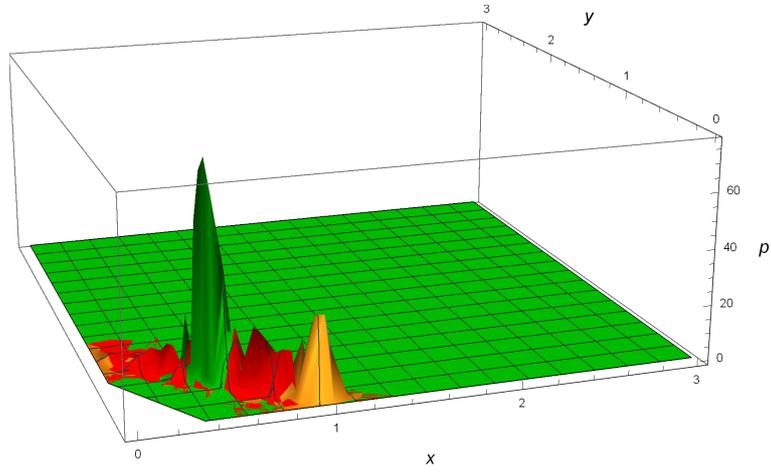}}}\hspace{4pt}
			\subfigure[FP equation for $t=1$ (Green), $t=2500$ (Red), $t=7000$ (Yellow), $\sigma_1=0.03$ and $\sigma_2=0.03$.]{
				\resizebox*{10cm}{!}{\includegraphics{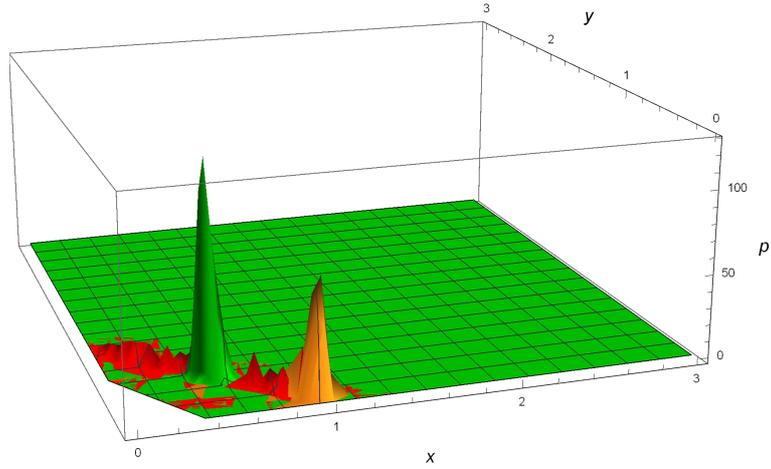}}}
			\caption{Numerical solution of the Fokker-Planck equation for $\theta=1$ and varying values of the noise intensity.}
		\end{center}
	\end{figure}
	
	The system was initialized such that $x,y$ have almost the same magnitude and such that the system starts on the line of steady states. In Figure 15 the evolution of the Fokker-Planck equation is seen for different times represented by different colours for the case of $\theta=1$ and varying noise intensities. In plot (b) we have the same noise intensities to compare with the following simulation of equation (40) in the next subsection. It is clear that for large time the system tends to $y=0$ and centred a bit off $x=1$.

	\begin{figure}[H]
		\begin{center}
			\subfigure[FP equation for $t=1000$, $\sigma_1=0.5$ and $\sigma_2=0.07$.]{
				\resizebox*{10cm}{!}{\includegraphics{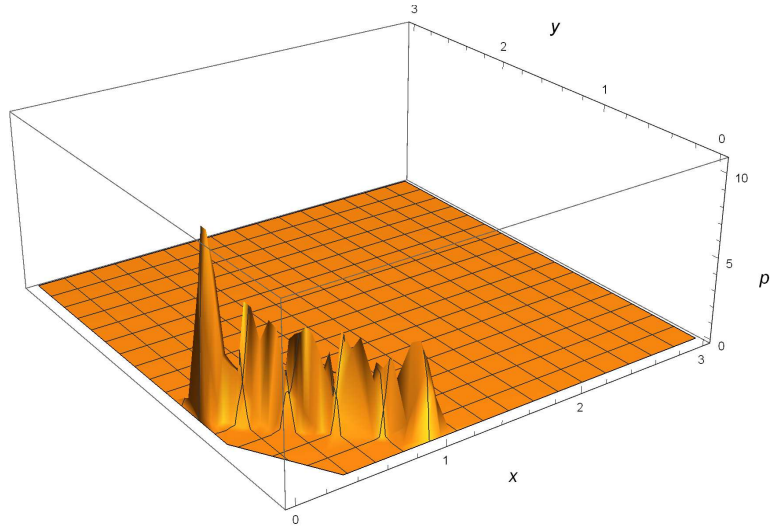}}}\hspace{4pt}
			\subfigure[FP equation for $t=1000$ and $\sigma_{1,2}=0.03$.]{
				\resizebox*{10cm}{!}{\includegraphics{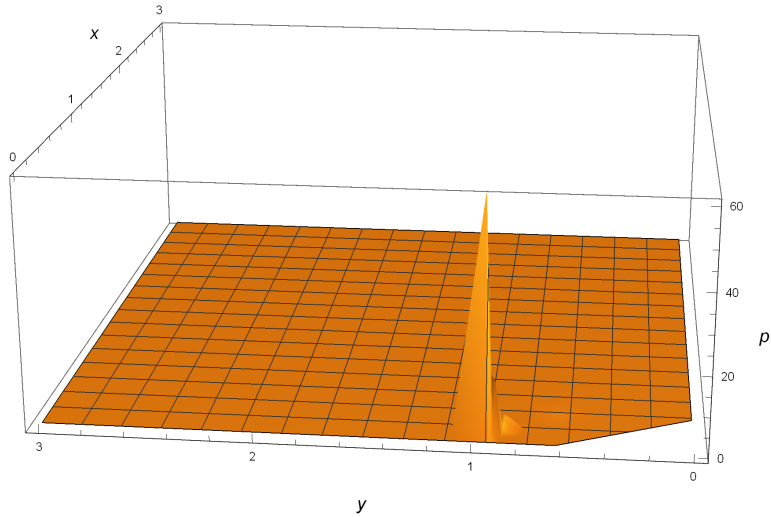}}}
			\caption{Numerical solution of the Fokker-Planck equation for $\theta=0.99$ and high and low noise intensity respectively.}
		\end{center}
	\end{figure}
	
	Figure 16(a) shows what happen in the case where $\theta=0.99$. As in the SDE case the probability density function implies that the system needs a lot of time to settle to one steady state. It seems that y population has a slight advantage here as we can also see in the numerical simulations of the Langevin equations before. Unfortunately due to numerical errors the solution of the Fokker-Planck dissipates so we are not able to properly explore the very long time behaviour of the solution which should probably be a peak around $y=1$. On the other hand for lower noise we get exactly what we would expect deterministically for the case of $\theta=0.99$ which is only for $y$ to survive (Fig 16b) which is the same as the result shown in the previous subsection.

	\subsection{Simulation of Differential-Algebraic system for dilution rate induced noise}
	
	As we already mentioned we need a bridge with the results from \cite{stephanopoulos_stochastic_1979}. In that paper the authors explored what happened in the case where z is simply assumed constant. In the present work, using the asymptotic analysis and the algebraic equation for z we made a simplification that allows us to see what happens when the variation of z is very small as well as the noise intensity low for the $\theta = 1$ scenario. It is easy to show that for the same level of noise intensity the variation of z in the case of (33-35) and (10-12) is very different. In the former the variation is significantly smaller as evident from the figure below.
	
	\begin{figure}[H]
		\begin{center}
			\subfigure[$\theta=1,\sigma=0.0006$.]{
				\resizebox*{10cm}{!}{\includegraphics{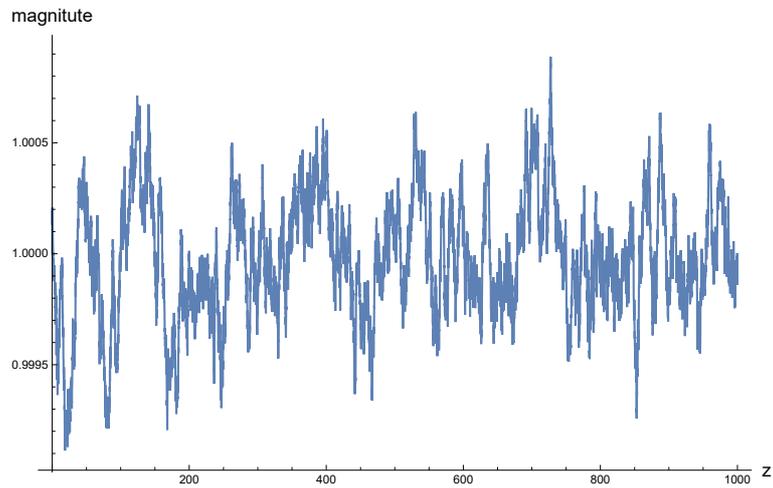}}}\hspace{4pt}
			\subfigure[$\theta=1,\sigma=0.0002$.]{
				\resizebox*{10cm}{!}{\includegraphics{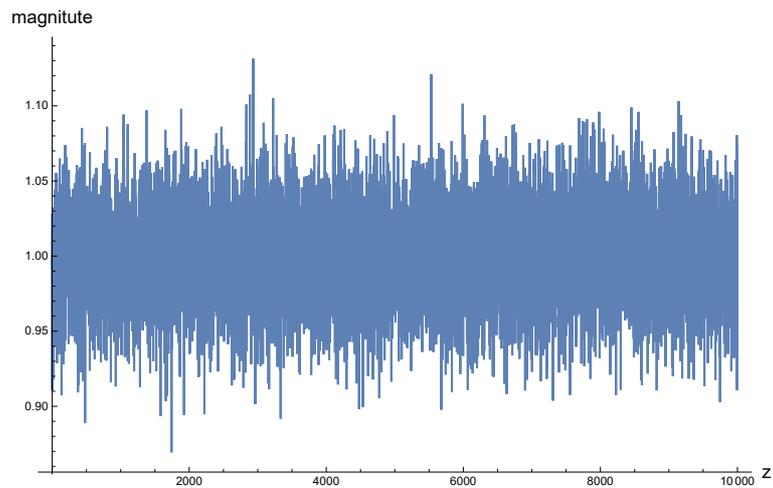}}}
			\caption{Evolution of z using the algebraic equation without noise and the differential equation with noise respectively. Notice how the fluctuation of the non-dimensional substrate are much three orders of magnitude smaller in the case where the algebraic equation is used.}
		\end{center}
	\end{figure}
	
	Again the system was initialized such that $x,y$ have almost the same magnitude and such that the system starts on the line of steady states. In Figure 18(a) the evolution of the Fokker-Planck equation is seen for different times represented by different colours for the case of $\theta=1$. It is clear that for large time the system tends to $y=0$ centred a bit off $x=1$. Figure 18(b) shows the large time behaviour of the system for $\theta=1$. Both plots are for low noise intensity with value $\sigma=0.03$.
	
	\begin{figure}[H]
		\begin{center}
			\subfigure[FP equation for $t=1$ (Green), $t=2500$ (Red) and $t=7000$ (Yellow).]{
				\resizebox*{10cm}{!}{\includegraphics{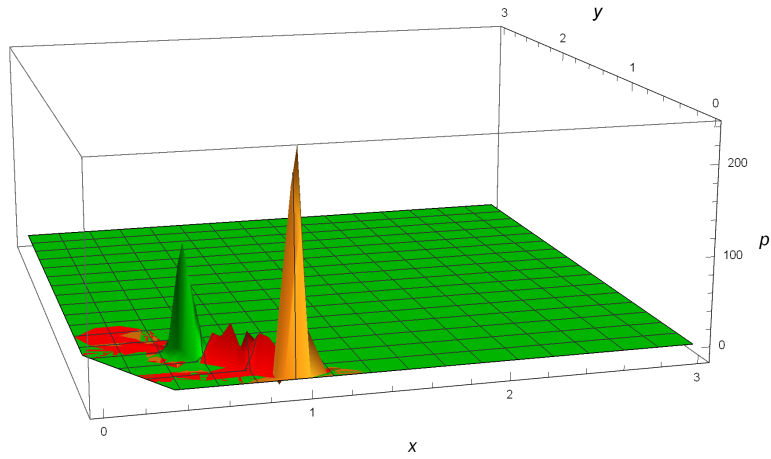}}}\hspace{4pt}
			\subfigure[FP equation for $t=7000$.]{
				\resizebox*{10cm}{!}{\includegraphics{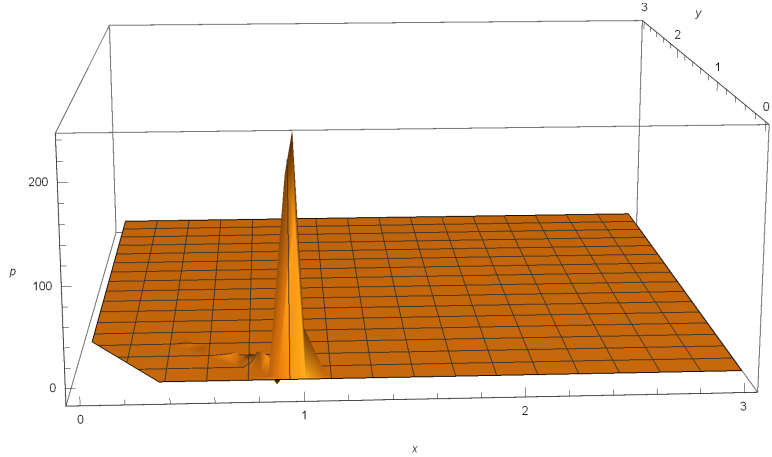}}}
			\caption{Numerical solution of the Fokker-Planck equation for $\theta=1$ and $\sigma=0.03$.}
		\end{center}
	\end{figure}
	
	Figure 19(a) demonstrated what happen in the case where $\theta=0.99$. We can see that the result is exactly what we expect for the case of low noise, i.e. $y$ reaches a steady-state which is again slightly smaller than 1 and $x$ dies out. On the other hand 19(b) shows that in the presence of high noise $\sigma=0.2$ this result changes and the population with the highest maximum growth rate is favoured.
	
	\begin{figure}[H]
		\begin{center}
			\subfigure[FP equation for $t=100$ and $\sigma=0.03$.]{
				\resizebox*{10cm}{!}{\includegraphics{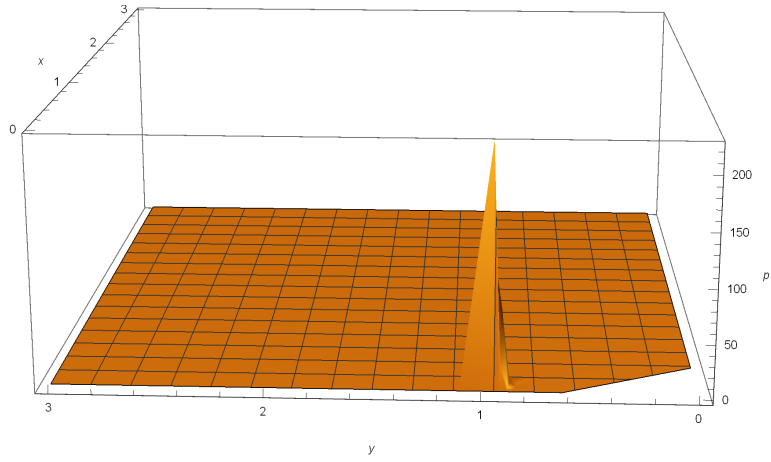}}}\hspace{4pt}
			\subfigure[FP equation for $t=300$ and $\sigma=0.2$.]{
				\resizebox*{10cm}{!}{\includegraphics{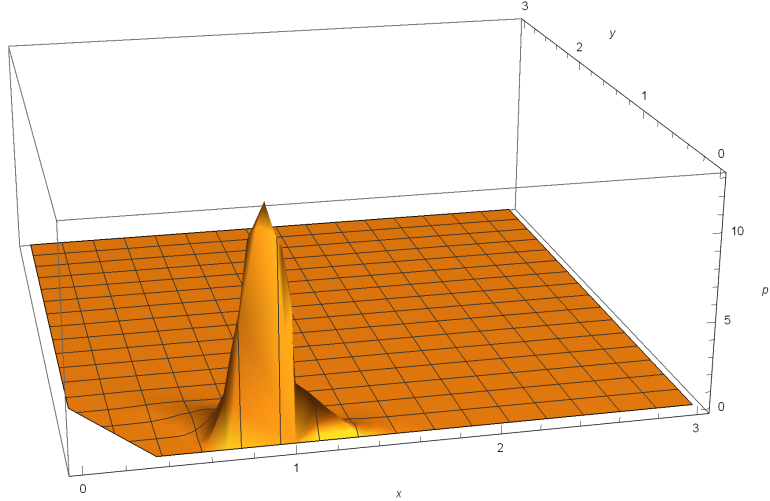}}}
			\caption{Numerical solution of the Fokker-Planck equation for $\theta=0.99$ and low and high noise intensity respectively.}
		\end{center}
	\end{figure}

	\subsection{Simulation of Langevin equations for full system with death}
	
	To explore the case where death of the microbial populations is included we picked some dummy values for the death rates such that there is a point of intersection between the modified growth curves of the two populations and such that the dominance for high values of z reverses. Meaning that as $z \to \infty$ $y$ has a greater growth rate than $x$ which is given by $a_2-\gamma_2$. In this case we have $a_2-\gamma_2 > a_1-\gamma_1$ and the relation between the parameter of the growth rate is given by (67). The simulation values chosen are:
	
	\begin{table}[ht]
		\centering
		\caption{Parameters and their values for the case of death rate.}
		{\begin{tabular}[l]{@{}lcccccc}\toprule
				Parameter (dimensionless) & Definition & Value \\ \hline
				$z_f$ & substrate feed & 15000 \\ 
				$a_i$ & maximum growth rate & 2.512, 1.411 \\ 
				$b_i$ & Michaelis constant & 0.041, 0.204 \\ 
				$\gamma_i$ & death rate & 1.41306, 0.171927\\ \hline
		\end{tabular}}
		\label{Table2}
	\end{table}
	
	The following figures are for both general noise and dilution rate induced noise respectively.
	
	\begin{figure}[H]
		\begin{center}
			\subfigure[$\theta=1,\sigma_1=0.0006$, $\sigma_2=0.0007$ and $\sigma_3=9$.]{
				\resizebox*{6cm}{!}{\includegraphics{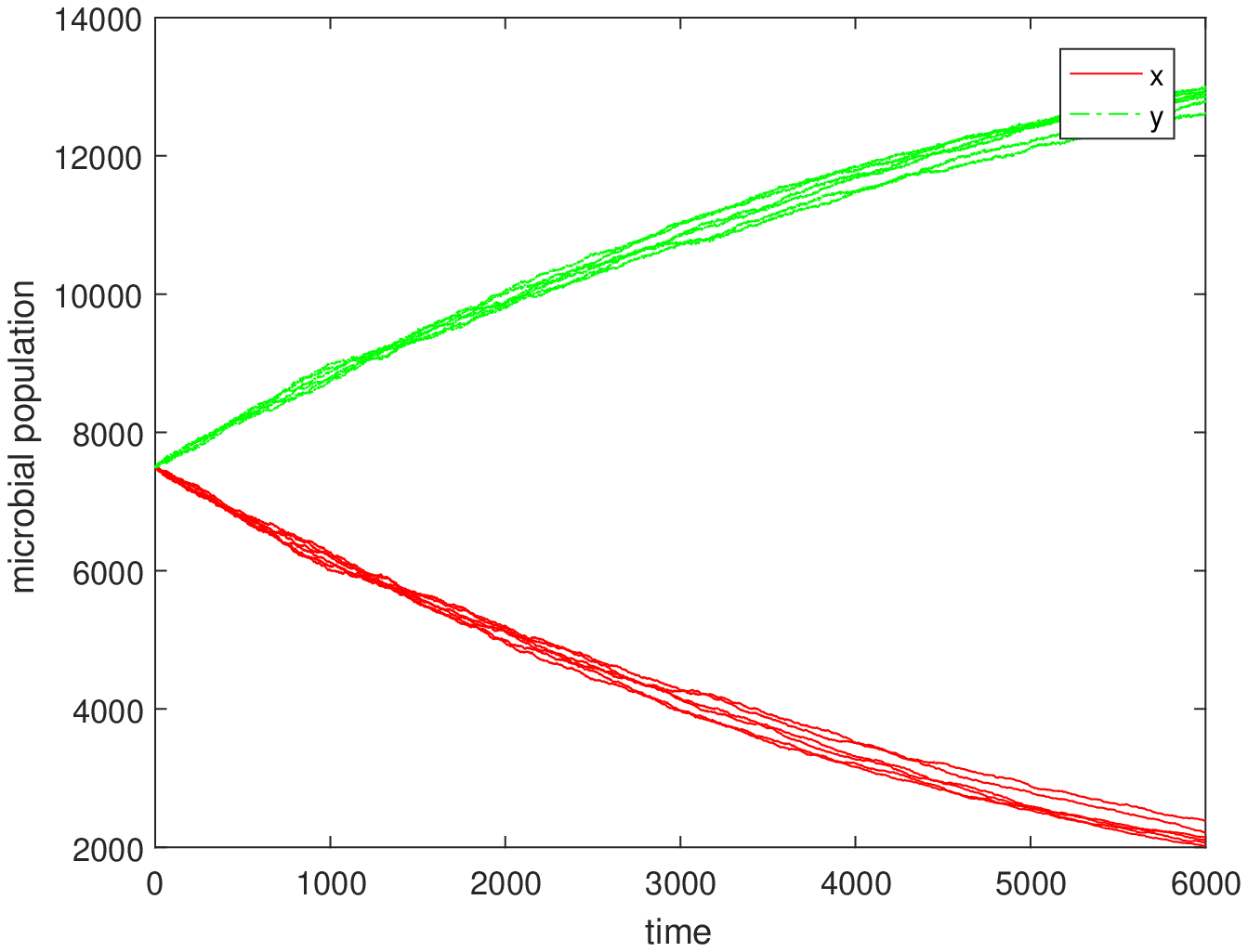}}}\hspace{5pt}
			\subfigure[$\theta=1,\sigma_1=0.05$, $\sigma_2=0.07$ and $\sigma_3=9$.]{
				\resizebox*{6cm}{!}{\includegraphics{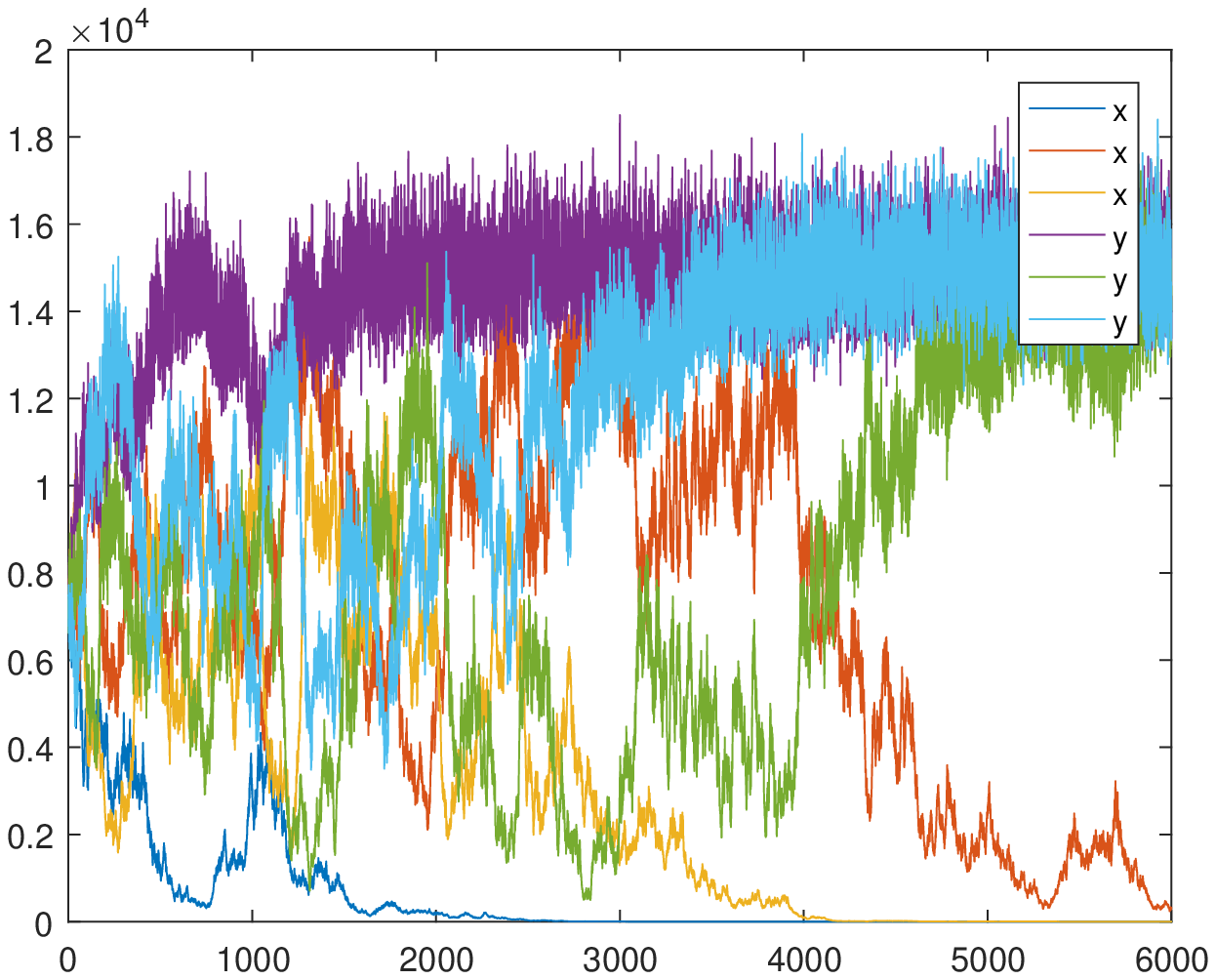}}}
			\caption{Plot of 6 runs of non-dimensional microbial populations versus non-dimensional time for the full system with death rate. and general noise}
		\end{center}
	\end{figure}
	
	\begin{figure}[H]
		\begin{center}
			\subfigure[$\theta=1,\sigma=0.0006$.]{
				\resizebox*{6cm}{!}{\includegraphics{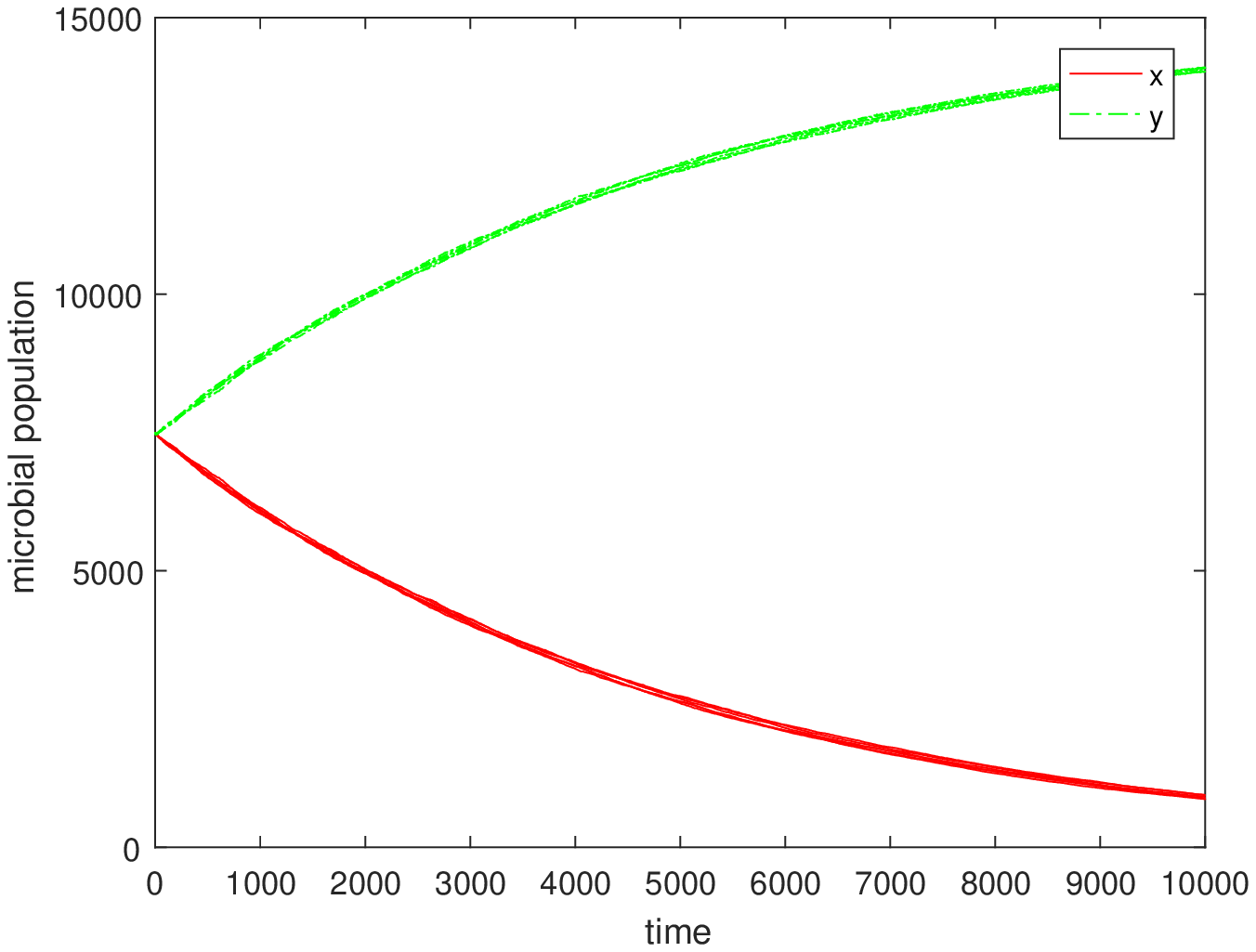}}}\hspace{5pt}
			\subfigure[$\theta=0.98,\sigma=0.03$.]{
				\resizebox*{6cm}{!}{\includegraphics{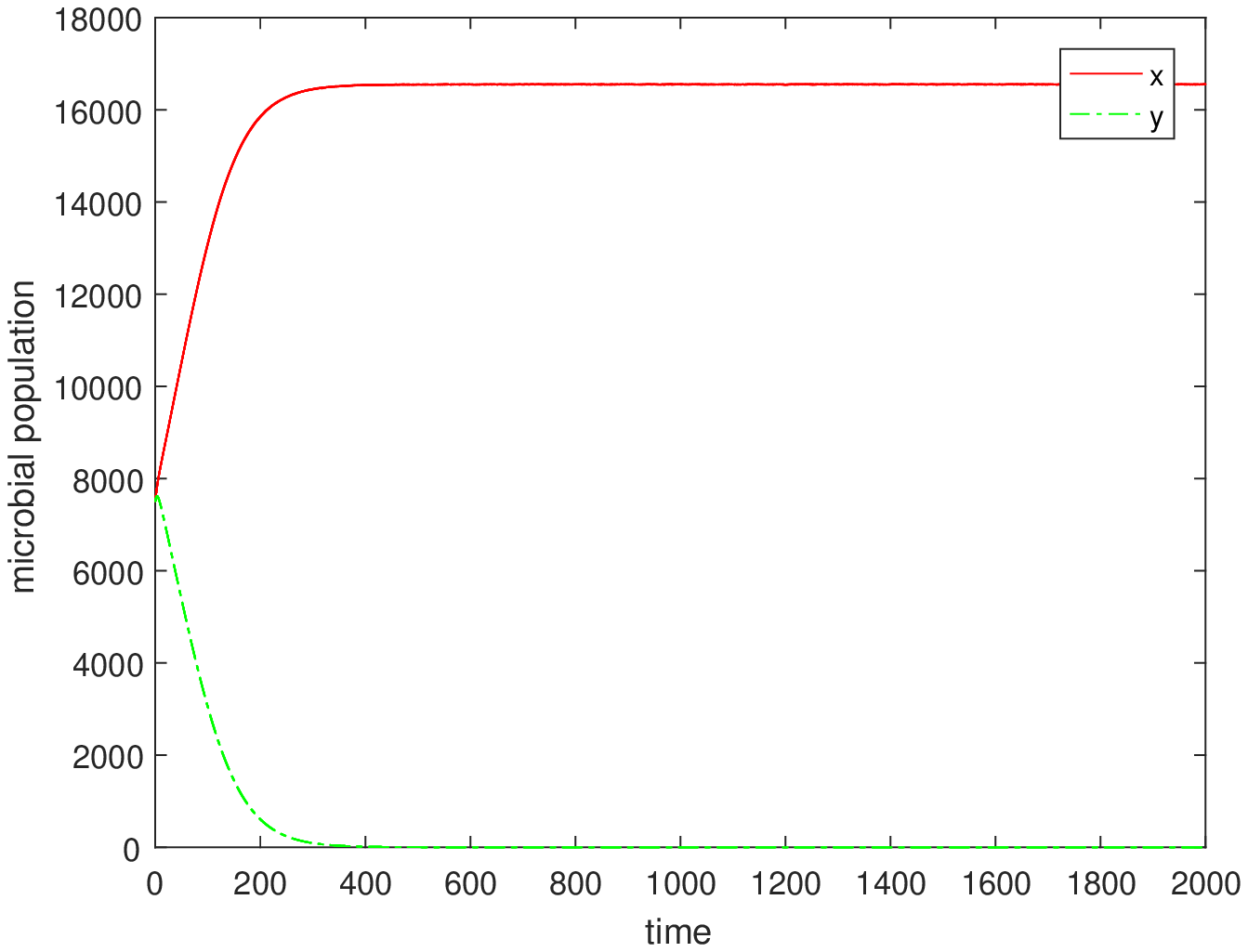}}}
			\caption{Plot of 6 runs of non-dimensional microbial populations versus non-dimensional time for the full system with death rate and dilution rate induced noise.}
		\end{center}
	\end{figure}
	
	Figures 20, 21 imply that when death rate is included the quantity that plays the most important role changes for the maximum growth rate ($a_i$) to that minus the death rate ($a_i-\gamma_i$). Other than that the results seem to be qualitative very similar to the respective noise cases without death. 
	
	\section{Discussion}
	
	In this project we investigated the coexistence of two populations in a chemostat competing over a single substrate under two distinct noise cases. the first was general environmental noise and the second was noise due to the fact that the chemostat dilution rate cannot be kept perfectly at a constant value. Extending the classical deterministic equations for the chemostat we included noise and formulated Ito type stochastic differential equations that were simulated at the deterministic steady-state line were the two populations are supposed to coexist. In both noise cases we found no evidence of coexistence and as always one population survived. In the dilution rate case the population with the higher maximum growth rate was always benefited independently of the noise intensity and in addition higher intensity seemed to help that population even further. The most important result was that for the dilution rate value where deterministically both populations would survive ($\theta=1$) only one did and it was always the same one. Similar results were found in the general noise case. No case where either population can survive was found as in \cite{stephanopoulos_stochastic_1979} although there might be some combination of $\sigma_{1,2}$ in the general noise scenario when that could possibly be observed.\\
	
	Following the results of the SDEs simulation a asymptotic analysis was performed which is valid when the dimensionless substrate feed is large. That asympotic analysis aided us in breaking the dynamics of the system in 5 stages where the last stage, being the one of interest, represents the approach of the system to the steady state and reduces the dimensions of our system to 2 instead of 3. That allows the formulation of a 2D Fokker-Planck equation which can be solved numerically. The resulting Fokker-Planck equation is slightly different depending on the noise assumptions and in both cases the results agree with the simulations of the SDEs. The benefit of formulating and Fokker-Planck equation is that if it is solved it admits the fate of the system without the need for multiple simulations and can clearly show the path of the system towards the final steady state.\\
	
	Finally, the effect of adding death rate was explored by first performing a linear stability analysis of the deterministic system in order to show that it remains that similar to the case without death rate. Following the stability analysis we performed simulations of the SDEs with dummy death rates such that the balance of the growth curves is reversed. Meaning that the for large substrate values the $y$ population had the highest growth rate which is given by the maximum growth rate minus the death rate as $z \to \infty$. Our results shows no qualitative difference with the no death rate cases which was expected.
	
	\section{Conclusion}
	
	Our work has been an investigation of coexistence and competitive exclusion in the case of general, non-linear Monod growth rates for two populations in a chemostat with stochastic effects taken into account. Based on the analysis and numerical simulation of two microbial populations competing in a chemostat in \cite{stephanopoulos_stochastic_1979} we wished to extend our understanding in two ways. The first was exploring the differences of solving the full system rather than a one-dimensional simplification and the second was adding a more general type of noise. In the original paper the results for $\theta=1$ are independent of noise intensity and suggest that either population has a probability of surviving depending solely on the initial conditions on the line of steady state.  Interestingly keeping the nature of the noise the same as the original paper but solving for the full system admits different results from the ones found in the one-dimensional case but only when $\theta=1$. That seems to remain true when we have general environmental noise.\\
	
	Despite the fact that in the case of general noise and high intensity the competition was more persistent something similar was not observed for low noise. Contrary to previous suggestions, one particular population becomes dominant at large times. An important step was the derivation of the Fokker-Planck equation for a reduced 2D system which can be used for more systematic exploration and simulation of it as well as the inclusion of death rate. Of course there is plenty of room for future work and further improvements. Currently our results are based on simulations of the respective stochastic systems so one possible path will be to try and prove the existence of the stochastic steady states analytically as was done for simpler linear growth rates \cite{xu_competition_2016} and potentially prove mathematically that in the case were $\theta=1$ there exist one or two steady states ($y=0$ or $x=0$) depending on the nature and intensity of the noise as well as the growth parameters. Additionally, a lot of work can be done regarding the Fokker-Planck equation derived here. As mentioned before due to numerical errors some probability density flux outwards of the domain leading to the total probability decreasing instead of staying constant at 1. A more systematic numerical solution of the equations could prevent that which would be useful in clarifying their long time behaviour and whether there is a bimodal distribution in the general noise case with $\theta=0.99$. Finally, one could explore the case were death rate is included and there are two intersection points in the upper right quadrant. It would be interesting to see how the fate of the system is affected when the two points are close and the noise is such that the dilution rate could potentially reach either value of the these two points.\\

\bibliographystyle{elsarticle-num} 
\bibliography{chemostatref}

	\begin{appendices}
		
		\section{Derivation of the Fokker-Planck equation from the Langevin SDEs}
		
		Without loss of generality let us take the one-dimensional case and assume that the evolution of the variable $f(t)$ is due to a deterministic, $G(f,t)$, and a stochastic element with multiplicative Gaussian white noise, $S(f,t)\zeta(t)$, such that $\e{\zeta(t_1)\zeta(t_2)} = \sigma^2\delta(t_1-t_2)$.
		
		\begin{equation}
		\frac{df(t)}{dt}=G(f,t) + S(f,t)\zeta(t)
		\end{equation}
		
		\noindent The above equation can be written in the form:
		
		\begin{equation}
		df(t)=G(f,t)dt + S(f,t)\sigma dW(t)
		\end{equation}
		
		\noindent Where, $dW(t)$ is a Wiener process increment.\\
		
		To formulate the corresponding Fokker-Planck equation we need to find the drift and diffusion terms, $\mu,D$ respectively. In order to achieve that we have to compute the moments of variable $f$ using the Langevin equation \cite{risken1984fokker}. The general formulas for the moments are:
		
		\begin{align}
		m^{(n)}(f',t)=&\frac{1}{n!}\lim_{dt \to 0}\frac{1}{dt}\int (f-f')^np(f,t+dt|f',t)\\
		&\frac{1}{n!}\lim_{dt \to 0}\frac{1}{dt}\e{(f(t+dt)-f(t))^n|f(t)=f'}
		\end{align}
		
		Using the fact that in (B.4) $f(t+dt)-f(t)=df(t)$ is given by (B.2) as well as the properties of the expectation of the Wiener process increments, namely $\e{dW(t)}=0,\e{(dw(t))^2}=dt$, we get:
		
		\begin{equation*}
		\e{df(t)|f(t)=f'}=\e{G(f,t)dt + S(f,t)\sigma dW(t)|f(t)=f'}=G(f',t)dt
		\end{equation*}
		\small
		\begin{align*}
		\e{(df(t))^2|f(t)=f'}\\
		=&\e{G^2(f,t)dt^2 + S^2(f,t)\sigma^2 dW^2(t) + 2G(f,t)S(f,t)\sigma dW(t)dt|f(t)=f'}\\
		=& S^2(f',t)\sigma^2dt +O((dt)^2)
		\end{align*}
		\normalsize
		\noindent Now, the drift and diffusion terms can be found:
		
		\begin{align*}
		\mu(f',t)=&\lim_{dt \to 0}\frac{1}{dt}\int (f-f')p(f,t+dt|f',t)\\
		=&\lim_{dt \to 0}\frac{1}{dt}\e{(f(t+dt)-f(t))|f(t)=f'} = G(f',t)
		\end{align*}
		
		\begin{align*}
		D(f',t)=&\lim_{dt \to 0}\frac{1}{2dt}\int (f-f')^2p(f,t+dt|f',t)\\
		=&\lim_{dt \to 0}\frac{1}{2dt}\e{(f(t+dt)-f(t))^2|f(t)=f'}=S^2(f',t)\sigma^2
		\end{align*}
		
		\noindent Hence, the corresponding Fokker-Planck is:
		
		\begin{equation*}
		\frac{\partial}{\partial t}p(f,t|f_0,t_0)=-\frac{\partial}{\partial f}\big[G(f,t)p(f,t|f_0,t_0)\big]+\frac{\sigma^2}{2}\frac{\partial^2}{\partial f^2}\big[S^2(f,t)p(f,t|f_0,t_0)\big]
		\end{equation*}
		\\
		
		The general formula for the n-dimensional Fokker-Planck equation is given by:
		\small
		\begin{equation}
		\frac{\partial}{\partial t}p(\bm{f},t|\bm{f_0},t_0)=-\sum_i\frac{\partial}{\partial f_i}\big[\mu(\bm{f},t)p(\bm{f},t|\bm{f_0},t_0)\big]+\sum_{i,j}\frac{\partial^2}{\partial f_if_j}\big[D_{i,j}(\bm{f},t)p(\bm{f},t|\bm{f_0},t_0)\big]
		\end{equation}
		\normalsize
		It can been shown \cite{ottinger_stochastic_1996} that for a general N-dimensional $It\bar{o}$ SDE system of the form:
		
		\begin{equation*}
		d\bm{X}_t=\bm{\mu}(\bm{X}_t,t)dt + \bm{\sigma}(\bm{X}_t,t)d\bm{W}_t
		\end{equation*}
		
		\noindent Where, $\bm{X}_t$ and $\bm{\mu}(\bm{X}_t,t)$ are N-dimensional vectors, $\bm{\sigma}(\bm{X}_t,t)$ is an $N$x$M$ matrix and $\bm{W}_t$ is an M-dimensional standard Wiener process.\\
		
		The corresponding Fokker-Planck of the above system is:
		\small
		\begin{equation}
		\frac{\partial p(\bm{x},t)}{\partial t}=-\sum_{i=1}^N\frac{\partial}{\partial x_i}\big[\mu_i(\bm{x},t)p(\bm{x},t)\big]+\frac{1}{2}\sum_{j=1}^M\sum_{i=1}^N\frac{\partial^2}{\partial x_i \partial x_j}\big[D_{i,j}(\bm{x},t)p(\bm{x},t)\big]
		\end{equation}
		\normalsize
		\noindent Here, $D_{i,j}=\sum_{k=1}^M\sigma_{ik}\sigma_{jk}$ or equivalently $\bm{D}(\bm{x},t)=\bm{\sigma}\cdot\bm{\sigma}^T $.
		
		\subsection{Application to chemostat model}
		
		The SDEs governing the stage 5 simplification of the chemostat model with general noise are given by:
		
		\begin{align*}
		d\bar{x}&=\bar{x}(f(\bar{x},\bar{y})-\theta_0)dt + \sigma_1 \bar{x} dW_1(t) ,\\
		d\bar{y}&=\bar{y}(g(\bar{x},\bar{y})-\theta_0)dt + \sigma_2 \bar{y} dW_2(t).\\
		\end{align*}
		
		\noindent Where, $N=M=2$ and
		
		\begin{equation*}
		\bm{\mu}(\bm{x},t)=\big(\bar{x}(f(\bar{x},\bar{y})-\theta_0),\bar{y}(g(\bar{x},\bar{y})-\theta_0)\big),
		\end{equation*}
		
		\[\sigma_{ij}= \left( \begin{array}{ccc}
		\sigma_1 \bar{x} & 0 \\
		0 & \sigma_2 \bar{y}  \\ \end{array} \right),\]\\
		
		\begin{equation*}
		\bm{W}_t=\big(W_1(t),W_2(t)\big).
		\end{equation*}
		
		Using the above it is trivial to find that the Fokker-Planck equation for our 2D system is:	
		
		\small
		
		\begin{align*}
		\frac{\partial}{\partial t}p(\bar{x},\bar{y},t|x_0,\bar{y}_0,t_0)=&-\frac{\partial}{\partial \bar{x}}\big[\bar{x}(f(\bar{x},\bar{y})-\theta_0)p(\bar{x},\bar{y},t|\bar{x}_0,\bar{y}_0,t_0)\big]\\
		&-\frac{\partial}{\partial \bar{y}}\big[\bar{y}(g(\bar{x},\bar{y})-\theta_0)p(\bar{x},\bar{y},t|\bar{x}_0,\bar{y}_0,t_0)\big]\\
		&+\frac{\sigma_1^2}{2}\frac{\partial^2}{\partial \bar{x}^2}\big[\bar{x}^2p(\bar{x},\bar{y},t|\bar{x}_0,\bar{y}_0,t_0)\big]+\frac{\sigma_2^2}{2}\frac{\partial^2}{\partial \bar{y}^2}\big[\bar{y}^2p(\bar{x},\bar{y},t|\bar{x}_0,\bar{y}_0,t_0)\big]
		\end{align*}
		
		\normalsize
		
		When noise originates from the dilution rate the increments and intensity are the same so we instead have:

		\[\sigma_{ij}= \left( \begin{array}{ccc}
		-\sigma \bar{x} & 0 \\
		-\sigma \bar{y} & 0  \\ \end{array} \right),\]\\
		
		\begin{equation*}
		\bm{W}_t=\big(W(t),W(t)\big).
		\end{equation*}
		
		From which we gain an extra term by using the methodology explained above:	
		
		\small
		
		\begin{align*}
		\frac{\partial}{\partial t}p(\bar{x},\bar{y},t|x_0,\bar{y}_0,t_0)=&-\frac{\partial}{\partial \bar{x}}\big[\bar{x}(f(\bar{x},\bar{y})-\theta_0)p(\bar{x},\bar{y},t|\bar{x}_0,\bar{y}_0,t_0)\big]\\
		&-\frac{\partial}{\partial \bar{y}}\big[\bar{y}(g(\bar{x},\bar{y})-\theta_0)p(\bar{x},\bar{y},t|\bar{x}_0,\bar{y}_0,t_0)\big]\\
		&+\frac{\sigma^2}{2}\frac{\partial^2}{\partial \bar{x}^2}\big[\bar{x}^2p(\bar{x},\bar{y},t|\bar{x}_0,\bar{y}_0,t_0)\big]+\frac{\sigma^2}{2}\frac{\partial^2}{\partial \bar{y}^2}\big[\bar{y}^2p(\bar{x},\bar{y},t|\bar{x}_0,\bar{y}_0,t_0)\big]\\
		&+\sigma^2\frac{\partial^2}{\partial \bar{x}\partial \bar{y}}\big[\bar{x}\bar{y}p(\bar{x},\bar{y},t|\bar{x}_0,\bar{y}_0,t_0)\big]
		\end{align*}

		\section{Equations, nondimensionalisation and stability analusis of system with death rate}
		The system of coupled equations describing the evolution of the microbial populations and the concentration of the substrate is the following:
		
		\begin{align}
		\dot{x_1}&=-\frac{qx_1}{V}+\mu_1(s)x_1,\\
		\dot{x_2}&=-\frac{qx_2}{V}+\mu_2(s)x_2,\\
		\dot{s}&=\frac{q(s_f-s)}{V}-\frac{1}{Y_1}\mu_1(s)x_1-\frac{1}{Y_2}\mu_2(s)x_2.
		\end{align}
		
		It is more convenient if equations (3.68-3.70) are nondimensionalised. Since the exploration of the intersection points is our main aim, the characteristic measures used for the nondimensionalisation will be the growth rate and the substrate values at the point of intersection. The values for which the growth rate curves intersect is $s_c$ so that $\mu_1(s_c)=\mu_2(s_c)=\mu_c$. The dimensionless variables are give below:
		
		\begin{equation*}
		x=\frac{x_1}{Y_1s_c}, y=\frac{x_2}{Y_2s_c}, z=\frac{s}{s_c}, z_f=\frac{s_f}{s_c}, \tau=\mu_ct, \theta=\frac{q}{V\mu_c}
		\end{equation*}
		
		\begin{equation*}
		f(z)=\frac{\mu_1(zs_c)}{\mu_c}, g(z)= \frac{\mu_2(zs_c)}{\mu_c}
		\end{equation*}
		
		By the definition of $z$ and the dimensionless growth rates, at the intersection point $z=1$ and $f(1)=g(1)=1$. The equations for the dimensionless growth rates are:
		
		\begin{align}
		f(z)=\frac{a_1z}{b_1+z}-\gamma_1\\
		g(z)=\frac{a_2z}{b_2+z}-\gamma_2
		\end{align}
		
		\noindent Where, $a_i=\frac{\mu_m^i}{mu_c}, b_i=\frac{K_s^i}{s_c}, \gamma_i=\frac{d_i}{\mu_c}$. Furthermore, due to the relation at the intersection point we obtain:
		
		\begin{equation}
		\gamma_i=\frac{a_i}{b_i+1}-1
		\end{equation}
		
		Finally, the new set of equations is:
		
		\begin{align}
		\dot{x}&=x(f(z)-\theta),\\
		\dot{y}&=y(g(z)-\theta),\\
		\dot{z}&=\theta(z_f-z)-xf(z)-yg(z).
		\end{align}
		
		\subsection{Linear stability analysis}
		Before proceeding to the stability analysis, the dimensions of the system can be reduced. Adding (68-70) admits:
		
		\begin{equation*}
		\dot{x}+\dot{y}+\dot{z}=\theta(z_f-z-x-y)
		\end{equation*}
		
		In steady state,
		
		\begin{equation}
		z=z_f-x-y
		\end{equation}
		
		\noindent This steady state is reached if the system is initialized inside the triangular domain defined by (71) \cite{stephanopoulos_stochastic_1979}. Now there is a new system of equations where (70) is replaced by (71). That system admits four steady states, (1) x=y=0 where both populations wash out, (2) x=0 and $g(z)=\theta$ where $x_1$ washes out and $x_2$ grows, (3) y=0 and $f(z)=\theta$ where $x_2$ washed out and $x_1$ grows and finally (4) $f(z)=g(z)=\theta$ where both populations coexist. The conditions for stability of these states will be explored. All the conditions found in the first section hold for the new growth parameters exactly as before.\\
		
		Two more conditions need to be taken into account before the stability analysis that those are the following:
		
		\begin{equation}
		\gamma_2>\gamma_1 \implies \frac{a_2}{b_2+1}>\frac{a_1}{b_1+1}
		\end{equation}
		
		\begin{equation}
		\gamma_i>0 \implies \frac{a_i}{b_i+1}>1
		\end{equation}
		
		\subsubsection{$x=y=0$}
		
		Linearising around the (0,0) steady-state gives the Jacobian matrix which can be used to find the eigenvalues. These are:
		
		\begin{equation*}
		\left\{\frac{a_2 z_f-(b_2+z_f) (\gamma_2+\theta )}{b_2+z_f},\frac{a_1 z_f-(b_1+z_f) (\gamma_1+\theta )}{b_1+z_f}\right\}
		\end{equation*}
		\\
		
		It is straightforward to show that these eigenvalues are negative when,
		
		\begin{equation*}
		f(z_f),g(z_f)<\theta.
		\end{equation*}
		
		\subsubsection{$x=0, g(z)=\theta$}
		
		Using the same approach of linearisation around the steady state we can find the eigenvalues for the second steady-state. In this case we have $x=0$ and we need to solve $g(z)=\theta$ to find the value for $y$. Specifically $y=z_f+b_2\frac{\theta+\gamma_2}{\theta+\gamma_2-a_2}$. The eigenvalues are:
		
		\small
		
		\begin{equation*}
		\left\{-\frac{(a_2-\gamma_2-\theta ) (a_2 z_f-(b_2+z_f) (\gamma_2+\theta ))}{a_2 b_2},\frac{(\gamma_2+\theta ) (a_1 b_2+(b_1-b_2) (\gamma_1+\theta ))-a_2 b_1 (\gamma_1+\theta )}{a_2 b_1-(b_1-b_2) (\gamma_2+\theta )}\right\}
		\end{equation*}
		\normalsize
		\\
		For the first eigenvalue to be negative two opposing conditions are found, $g(z_f)>\theta$. The last means that $\theta$ is greater than the value of $g(z)$ as $z \rightarrow \infty$. \\
		
		The conditions for the second eigenvalue to be negative are $\theta<1$. \\
		
		\subsubsection{$y=0, f(z)=\theta$}
		
		Linearising around this steady state gives the eigenvalues,
		
		\small
		
		\begin{equation*}
		\left\{-\frac{(a_1-\gamma_1-\theta ) (a_1 z_f-(b_1+z_f) (\gamma_1+\theta ))}{a_1 b_1},\frac{(\gamma_2+\theta ) (b_2(-a_1+\gamma_1+\theta )-b_1 (\gamma_1+\theta ))+a_2 b_1 (\gamma_1+\theta )}{a_1 b_2+(b_1-b_2) (\gamma_1+\theta )}\right\}
		\end{equation*}
		\normalsize
		\\
		
		Condition for the first eigenvalue to be negative are $f(z_f)>\theta$. As for the second eigenvalue the condition is simply $\theta>1$.\\
		
		\subsubsection{$f(z)=g(z)=\theta$}
		
		In the final steady-state we have co-existence and the growth rates, the dimensionless dilution $(\theta)$ rate as well as $z$ are all equal to 1. Here, there is a line of fixed-points and not a single point. That line is given by the equation $y=z_f-x-1$. Since the intersection points exists and is at the positive quadrant the line must also be in the upper right $(x,y)$ quadrant. So, $z_f>1$, otherwise the line will never be in a region where both $x,y>0$ which is the biologically relevant region. Linearising the system around an arbitrary point (A,B) where $B=z_f-A-1$, admits these eigenvalues:
		
		\begin{equation*}
		\{\frac{2 a_1 b_1 (b_2+1)^2 (A-z_f+1)-2 A a_2 (b_1+1)^2 b_2}{2 (b_1+1)^2 (b_2+1)^2},0\}
		\end{equation*}
		\\
		Since A is part of the fixed-points line we know that $A<z_f-1$ and we can show that the first eigenvalue is always negative.
		
	\end{appendices}

\end{document}